\documentclass{amsart}
\let\<=\langle
\let\>=\rangle
\let\\=\backslash
\let\vphi=\varphi
\let\noi=\noindent
\let\sse=\subseteq
\let\veps=\varepsilon
\let\limply=\Longrightarrow

\def\0{\{0\}}
\def\span{{\kern.5pt{\rm span}\kern1pt}}
\def\smallfrac#1#2{{\textstyle{\frac{#1}{#2}}}}
\def\conv{{\;\to\;}}
\def\uconv{{\kern3pt{\buildrel_{\scriptstyle u}\over\longrightarrow}\kern3pt}}
\def\sconv{{\kern2pt{\buildrel_{\scriptstyle s}\over\longrightarrow}\kern3pt}}
\def\wconv{{\kern2pt{\buildrel_{\scriptstyle w}\over\longrightarrow}\kern3pt}}
\def\notconv{{{\kern-1pt\conv\kern-11pt\slash}\kern6pt}}

\def\smallmatrix#1{\null\,\vcenter{
                   \baselineskip=8pt\mathsurround=0pt\ialign{
                   \hfil ${\scriptstyle##}$
                   \hfil &&
                   \hfil ${\scriptstyle##}$
                   \hfil \crcr
                   \mathstrut \crcr
                   \noalign{\kern-\baselineskip}#1 \crcr
                   \mathstrut \crcr
                   \noalign{\kern-\baselineskip} \crcr }}\!}

\def\B{{\mathcal B}}
\def\M{{\mathcal M}}
\def\N{{\mathcal N}}
\def\R{{\mathcal R}}
\def\X{{\mathcal X}}
\def\Y{{\mathcal Y}}
\def\BX{{\B[\X]}}

\def\CC{{\mathbb C\kern.5pt}}
\def\RR{{\mathbb R\kern.5pt}}

\newtheorem{theorem}{Theorem}[section]
\newtheorem{corollary}{Corollary}[section]
\newtheorem{proposition}{Proposition}[section]
\newtheorem{lemma}{Lemma}[section]
\theoremstyle{definition}
\newtheorem{remark}{Remark}[section]

\begin{document}

\vglue-50pt\noi
\hfill{\it Advances in Mathematical Sciences and Applications}\/
{\bf 29} (2020) 145--170

\vglue30pt
\title
{Asymptotic Limits, Banach Limits, and Ces\`aro Means}
\author{C.S. Kubrusly}
\address{Mathematics Institute, Federal University of Rio de Janeiro, Brazil}
\email{carloskubrusly@gmail.com}
\author{B.P. Duggal}
\address{Faculty of Sciences and Mathematics, University of Ni\v s, Serbia}
\email{bpduggal@yahoo.co.uk}
\renewcommand{\keywordsname}{Keywords }
\date{April 20, 2020}

\begin{abstract}
Every new inner product in a Hilbert space is obtained from the original one
by means of a unique positive operator$.$ The first part of the paper is a
survey on applications of such a technique, including a characterization of
similarity to isometries$.$ The second part focuses on Banach limits for
dealing with power bounded operators. It is shown that if a power bounded
operator for which the sequence of shifted Ces\`aro means converges (at least
in the weak topology) uniformly in the shift parameter, then it has a Ces\`aro
asymptotic limit coinciding with its $\vphi$-asymptotic limit for all Banach
limits $\vphi$.
\end{abstract}

\maketitle

\section{Introduction}

The purpose of this paper is twofold$.$ It is a survey with an expository
flavor linking the notions in the title, and also includes original results.

\vskip6pt
The paper is split into two parts, both dealing with bounded linear operators
on a Hilbert space$.$ The first part (Sections 3, 4 and 5) surveys the
technique of gen\-erating a new inner product from the original one, and its
applications to sim\-ilarity to isometries and asymptotic limit for
contractions, emphasizing the common role played by the equation
${T^*\!A\,T\kern-1pt=\kern-1ptA}.$ The central results of this part appear in
Prop\-ositions 4.1, 4.2, giving a comprehensive characterization of similarity
to isometries.

\vskip6pt
The second part (Sections 6 and 7) is a follow-up of the first one, extending
it to power bounded operators by means of the $\vphi$-asymptotic limit
associated with a \hbox{Banach} limit $\vphi$ and, alternatively, to bounded
operators by means of Ces\`aro asymptotic limit associated with Ces\`aro
means, still focusing on the role played by the equation ${T^*\!A\,T=A}.$
Theorem 6.1 brings together a large collection of properties of
$\vphi$-asymptotic limits for power bounded operators, as a generalization of
analogous results for contractions$.$ Similarly, Theorem 7.1 brings together
a large collection of properties of Ces\`aro-asymptotic limits for bounded
operators$.$ Theorem 7.2 shows that if a power bounded operator is such that
its sequence of Ces\`aro means converges in the weak topology, whose shifted
sequences converge uniformly in the shift parameter, then its Ces\`aro
asymptotic limit coincides with its $\vphi$-asymptotic limit for all Banach
limits $\vphi.$ This is followed by an application in Corollary 7.1.

\section{Notation and Terminology}

A linear transformation $L$ on a linear space $\X$ is injective if and only if
its kernel ${\N(L)=L^{-1}(\0)}$ is null (i.e., if and only if ${\N(L)=\0}).$
If $\X$ is a normed space, then let $\BX$ stand for the normed algebra of all
operators on $\X$ (i.e., of all bounded linear transformations of $\X$ into
itself)$.$ If ${T\kern-1pt\in\kern-1pt\BX}$, then $\N(T)$ is a subspace of
$\X$, which means a {\it closed}\/ linear manifold of $\X.$ The range
${\R(T)=T(\X)}$ of ${T\kern-1pt\in\kern-1pt\BX}$ is a (not necessarily closed)
linear manifold of $\X.$ An operator $T$ on a normed space $\X$ has a bounded
inverse on its range if and only if it is bounded below$.$ An operator $T$ on
a Banach space $\X$ is bounded below if and only if it is injective with a
closed range (i.e., ${\N(T)=\0}$ and ${\R(T)=\R(T)^-}\!$ where the upper bar
denotes closure).

\vskip6pt
Let $\X$ be a normed space$.$ An operator ${T\kern-1pt\in\kern-1pt\BX}$ is an
isometry if ${\|Tx\|=\|x\|}$ for every ${x\kern-1pt\in\kern-1pt\X}$ (a unitary
operator is an invertible isometry on a Hilbert space)$.$ It is a contraction
if ${\|Tx\|\kern-1pt\le\kern-1pt\|x\|}$ for every ${x\kern-1pt\in\X}$ (i.e.,
${\|T\|\kern-1pt\le\kern-1pt1}$), and it is power bounded if
${\sup_n\|T^n\|\kern-1pt<\kern-1pt\infty}.$ In this case set
${\beta=\sup_n\|T^n\|}.$ Thus $T$ is power bounded if there is a constant
${\beta\kern-1pt>\kern-1pt0}$ such that ${\|T^nx\|\le\beta\|x\|}$ for all
integers ${n\kern-1pt\ge\kern-1pt1}$ and every ${x\kern-1pt\in\kern-1pt\X}$,
which implies ${\sup_n\|T^nx\|\kern-1pt<\kern-1pt\infty}$ for every
${x\kern-1pt\in\kern-1pt\X}.$ The converse holds if $\X$ is a Banach space by
the Banach--Steinhaus Theorem$.$ Every isometry is a contraction and every
contraction is power bounded$.$ An operator ${T\kern-1pt\in\kern-1pt\BX}$ is
power bounded below if there is a constant ${\alpha\kern-1pt>\kern-1pt0}$ such
that ${\alpha\|x\|\le\|T^nx\|}$ for all integers ${n\kern-1pt\ge\kern-1pt1}$
and every ${x\kern-1pt\in\kern-1pt\X}.$ A $\BX$-valued sequence $\{S_n\}$
converges uniformly (or in the operator norm topology) to an operator
${S\kern-1pt\in\kern-1pt\BX}$ if ${\|(S_n\!-S)\|\to0}$ (notation$:$
${S_n\kern-1pt\uconv S}).$ It converges strongly to $S$ if the $\X$-valued
sequence $\{S_nx\}$ converges to $Sx$ in the norm topology (i.e.,
${\|(S_n\!-S)x\|\to0}$) for every ${\kern-1ptx\in\kern-1pt\X}$ (notation$:$
${S_n\kern-1pt\sconv S}).$ The sequence $\{S_n\}$ converges weakly to
${S\kern-1pt\in\kern-1pt\BX}$ if ${f((S_n\!-S)x)\to0}$ for every $f$ in the
dual $\X^*\!$ of $\X$ and every $x$ in $\X$ (notation$:$
${S_n\kern-1pt\wconv S}$ --- if $\X$ is a Hilbert space with inner product
${\<\cdot\,;\cdot\>}$, weak convergence means ${\<(S_n\!-S)x\,;y\>\to0}$ for
every ${x,y\in\X}$ by the Riesz Representation Theorem, which is equivalent to
${\<(S_n\!-S)x\,;x\>\to0}$ for every ${x\in\X}$ if the Hilbert space is
complex by the Polarization Identity)$.$ Uniform convergence clearly implies
strong convergence, which in turn implies weak convergence$.$ An operator
${\!T\kern-1pt\in\kern-1pt\BX}$ is of class $C_{0{\textstyle\cdot}}$ if the
power sequence $\{T^n\}$ converges strongly to the null operator,
${T^nx\conv0}$ for every ${x\kern-1pt\in\kern-1pt\X}$ (i.e., if $T$ is
strongly stable --- notation$:$ ${T^n\sconv O}$), and it is of Class
$C_{1{\textstyle\cdot}}$ if ${T^nx\notconv0}$ for every
${0\kern-1pt\ne\kern-1ptx\kern-1pt\in\kern-1pt\X}$.

\vskip6pt
Suppose $\X$ is an inner product space with inner product
$\<\,\cdot\,;\cdot\,\>.$ The norm induced by the inner product will be denoted
by $\|\cdot\|.$ If $\X$ is a Hilbert space and ${T\kern-1pt\in\kern-1pt\BX}$,
then ${T^*\!\in\BX}$ denotes its (Hilbert-space) adjoint$.$ A self-adjoint
operator $A$ (i.e., one for which ${A^*\!=A}$) is nonnegative or positive if,
respectively, ${0\le\<Ax\,;x\>}$ for every ${x\in\X}$ or ${0<\<Ax\,;x\>}$ for
every nonzero ${x\in\X}$ (notation$:$ ${A\ge O}$ or ${A>O}).$ A (self-adjoint)
operator $A$ is positive if and only if it is nonnegative and injective$:$
$$
A>O
\iff
A\ge O
\;\;\hbox{and}\;\;
\N(A)=\0.
$$
Injective self-adjoint operators have dense range (i.e., $\R(A)^-\!=\X$
whenever $\N(A)=\0$ if $A^*\!=A).$ Thus positive operators are injective with
dense range$.$ Hence a positive operator $A$ on a Hilbert space $\X$ is
bounded below if and only if its has a bounded inverse on its closed dense
image, which in turn is equivalent to saying that it is injective and
surjective, which means invertible (with a bounded inverse)$.$ Invertible
positive operators are called strictly positive and denoted by ${A\succ O}$:
$$
A\succ 0
\iff
\hbox{$A>O$ has a bounded inverse on $\X$}.
$$
A nonnegative operator $A$ has a unique nonnegative square root
$A^\frac{1}{2}$ which is positive or strictly positive whenever $A$ is
(indeed, $\N(A^\frac{1}{2})=\N(A)$ and $\R(A^\frac{1}{2})^-\!=\R(A)^-)$.

\section{Generating a New Inner Product}

Take a linear space $\X$, let ${\<\,\cdot\,;\cdot\,\>}$ be an inner product
on $\X$, suppose ${(\X,\<\,\cdot\,;\cdot\,\>)}$ is a Hilbert space, and let
$A$ be a nonnegative operator (i.e., ${A\ge O}$) on this Hilbert space$.$ As
is readily verified, the nonnegative operator $A$ generates a new semi-inner
product ${\<\,\cdot\,;\cdot\,\>_A}$ on $\X$ defined for every ${x,y\in\X}$ by
$$
\<x\,;y\>_A=\<Ax\,;y\>,
$$
which induces a seminorm ${\|\cdot\|_A}$ on $\X$ given by
$$
\|x\|_A=\<Ax\,;x\>^\frac{1}{2}=\|A^\frac{1}{2}x\|.
$$
This seminorm $\|\cdot\|_A$ becomes a norm whenever $A$ is injective (i.e.,
${\N(A)=\0}$ or, equivalently, whenever the nonnegative $A$ is positive)$.$
Consider the inner product space ${(\X,\<\,\cdot\,;\cdot\,\>_A)}.$ As defined
in Section 2, $\BX$ is the Banach algebra of all linear operators on the
Hilbert space ${(\X,\<\,\cdot\,;\cdot\,\>)}$ which are bounded (i.e.,
continuous) with respect to the norm $\|\cdot\|$ induced by the inner product
${\<\,\cdot\,;\cdot\,\>}.$ If $A$ is positive, then let $\BX_A$ denote the
normed algebra of all linear operators on the inner product space
${(\X,\<\,\cdot\,;\cdot\,\>_A)}$ which are bounded (i.e., continuous) with
respect to the new norm $\|\cdot\|_A.$ In this case (i.e., if ${A>O}$ or,
equivalently, if $\|\cdot\|_A.$is a norm), the following elementary
result represents an appropriate starting point.

\vskip3pt\noi
\begin{proposition}
Let\/ ${A>O}$ be an arbitrary positive operator on a Hilbert space\/
${(\X,\<\,\cdot\,;\cdot\,\>)}$ and consider the norm\/ ${\|\cdot\|_A}$ induced
by the inner product\/ ${\<\,\cdot\,;\cdot\,\>_A}={\<A\,\cdot\,;\cdot\,\>}.$
The following assertions are equivalent.

\vskip4pt\noi
\begin{description}
\item{$\kern-9pt$\rm(a)$\kern1pt$}
$A$\/ is invertible\quad
{\rm(i.e.,} has bounded inverse on\/ $\X$; equivalently,\/ ${A\succ O})$.

\vskip6pt
\item{$\kern-9.5pt$\rm(b)$\kern0pt$}
The norms\/ ${\|\cdot\|_A}$ and\/ ${\|\cdot\|}$ on\/ $\X$ are equivalent.
\end{description}
\end{proposition}

\vskip4pt\noi
{\it Proof.}\/
Take the Hilbert space ${(\X,\<\,\cdot\,;\cdot\,\>)}$ and let $A$ be a
positive operator on it$.$ So
$$
\|Ax\|^2\!
\le\|A^{\frac{1}{2}}\|^2\|A^{\frac{1}{2}}x\|^2\!
=\|A\|\,\|x\|_A^2
=\|A\|\,|\<Ax\,;x\>|
\le\|A\|\,\|Ax\|\,\|x\|
\le\|A\|^2\|x\|^2\!
$$
for every $x$ in $\X.$ Since $A$ is positive, it is injective, and so $A$ is
bounded below if and only if it has a closed range$.$ Since $A$ is an
injective self-adjoint, then ${\R(A)^-\!=\X}.$ Thus the positive operator $A$
is invertible (i.e., ${A>O}$ has a bounded inverse on $\X$ or, equivalently,
$A$ is bounded below and surjective) if and only if $A$ is bounded below,
which means $\alpha^2\|x\|^2\le\|Ax\|^2$ for every ${x\in\X}$ and some
${\alpha>0}.$ Therefore
$$
\kern-10pt
A\succ O
\;\iff\;
\smallfrac{\alpha^2}{\|A\|}\|x\|^2\le\|x\|_A^2\le\|A\|\|x\|^2\!
\;\;\hbox{for every $x\in\X$, for some $\alpha>0$}. \eqno{\qed}
$$

\vskip6pt
Perhaps the first result along this line is one ensuring that for every new
inner product there is a positive operator generating it$.$ $\!$This is a
classical result from \cite{Sto}.

\vskip3pt\noi
\begin{lemma}
Let\/ ${(\X,\<\,\cdot\,;\cdot\,\>)}$ be a Hilbert space, and let\/
${[\,\cdot\,;\cdot\,]}$ be a semi-inner product in\/ $\X.$ Then there exists
a unique nonnegative operator\/ ${A\in\BX}$ for which
$$
[x\,;y]=\<x\,;y\>_A=\<Ax\,;y\>
\;\;\,\hbox{for every}\;\;
x,y\in\X.
$$
If this unique\/ ${A\in\BX}$ is positive, then\/ ${[\,\cdot\,;\cdot\,]}$
becomes an inner product in\/ $\X$.
\end{lemma}

\vskip0pt\noi
\begin{proof}
This is a particular case of a fundamental result for densely defined bounded
sesquilinear forms ${[\,\cdot\,;\cdot\,]}$ in a Hilbert-space setting
\cite[Theorem 2.28, p.63]{Sto}$.$ In particular, if
${(\X,\<\,\cdot\,;\cdot\,\>)}$ is a Hilbert space, then the result holds for
every Hermitian symmetric sesquilinear form inducing either a nonnegative or a
positive quadratic form (i.e., it holds for every semi-inner or inner product
${[\,\cdot\,;\cdot\,]}$ on ${\X\times\X}$) according to whether $A$ is
nonnegative or positive, respectively.
\end{proof}

\vskip6pt
The inner product space ${(\X,[\,\cdot\,;\cdot\,])}$ may not be a Hilbert
space if the positive $A$ is not strictly positive$.$ However, if $A$ is
invertible (i.e., if ${A\succ O}$), then the new norm generated by $A$ is
equivalent to the original one (by Proposition 3.1), and so
${(\X,\<\,\cdot\,;\cdot\,\>_A)}$ becomes a Hilbert space.

\section{Similarity to an Isometry and the Equation ${T^*\!A\,T=A}$}

Let the seminorm (norm) $\|\cdot\|_A$ be the one induced by the new semi-inner
product (inner product) ${\<\cdot\,;\cdot\>_A}={\<A\cdot\,;\cdot\>}$ as
discussed in the previous section.

\vskip3pt\noi
\begin{proposition}
Let\/ $\X$ be a Hilbert space$.$ Take an arbitrary operator\/ $T$ in\/ $\BX$
and an arbitrary nonnegative operator\/ $A$ in\/ $\BX$.

\vskip4pt\noi
\begin{description}
\item{$\kern-8pt$\rm(a)$\kern1pt$}
Since\/ ${A\ge O}$, then $\|Tx\|_A=\|x\|_A$ for every\/ ${x\in\X}$ if and only
if\/ ${T^*\!A\,T=A}$.

\vskip6pt
\item{$\kern-8pt$\rm(b)$\kern0pt$}
If\/ ${A>O}$, then\/ $T$ is an isometry in\/ $\BX_A$ if and only if\/
${T^*\!A\,T=A}$.

\vskip6pt
\item{$\kern-8pt$\rm(c)$\kern1pt$}
If\/ ${A\succ O}$ and ${T^*\!A\,T=A}$, then $T$ is similar to an isometry.

\vskip6pt
\item{$\kern-8pt$\rm(d)$\kern.5pt$}
If\/ $T$ is similar to an isometry, then\/ ${T^*\!A'\kern1ptT=A'}$ for
some\/ ${O\prec A'\in\BX}$.
\end{description}
\end{proposition}

\vskip0pt\noi
\begin{proof}
(a)
If ${A\ge O}$, then $\|\cdot\|_A$ is a semi-norm on $\X.$ Since
${T^*\!A\,T\kern-1pt-\kern-1ptA}$ is self-adjoint,
${\<(T^*\!A\,T\kern-1pt-\kern-1ptA)x,x\>}=0$ for every ${x\in\X}$ if and only
if ${T^*\!A\,T=A}.$ Therefore, since
$$
\|Tx\|_A^2\!=\!\|A^\frac{1}{2}Tx\|^2\!=\!\<A\,Tx\,;Tx\>\!=\!\<T^*\!A\,Tx,x\>
\quad\;\;\hbox{and}\;\;\quad
\<Ax,x\>\!=\!\|A^\frac{1}{2}x\|^2\!=\!\|x\|_A^2
$$
\vskip2pt\noi
for every ${x\in\X}$, we get the result in (a).

\vskip6pt\noi
(b)
If ${A>O}$, then $\|\cdot\|_A$ is a norm on $\X$ and therefore the identity
${\|Tx\|_A=\|x\|_A}$ for every ${x\in\X}$ means the operator $T$ is an
isometry in $\BX_A.$ Now apply (a).

\vskip6pt\noi
(c)
If ${A\ge O}$ and ${T^*\!A\,T=A}$, then for every ${x\in\X}$
$$
\|A^\frac{1}{2}Tx\|^2
=\<{A^\frac{1}{2}Tx\,;A^\frac{1}{2}Tx}\>
=\<{T^*\!A\,Tx\,;x}\>
=\<{Ax\,;x}\>
=\|A^\frac{1}{2}x\|^2.
$$
If ${A\succ O}$, then $\|A^\frac{1}{2}TA^{-\frac{1}{2}}x\|=\|x\|$ for every
${x\in\X}.$ So $T$ is similar to an isometry.

\vskip6pt\noi
(d)
If ${T\in\BX}$ is similar to an isometry, then there exists an invertible
transformation $W$ in ${\B[\X,\Y]}$ (with a bounded inverse $W^{-1}$ in
${\B[\Y,\X]}$ for some Hilbert space $\Y$ unitarily equivalent to $\X$ by the
Inverse Mapping Theorem since $\X$ is Banach) for which ${WTW^{-1}}$ is an
isometry in $\B[\Y].$ Since $W$ is invertible, the polar decomposition of $W$
is given by $W=U|W|$ where ${U\in\B[\X,\Y]}$ is unitary and
$|W|={(W^*W)^\frac{1}{2}\in\BX}$ is strictly positive$.$ Thus
${O\prec|W|=U^*W}.$ Since ${WTW^{-1}}$ is an isometry on $\Y$, then
${U^*WTW^{-1}U}$ is an isometry on $\X.$ Thus
\begin{eqnarray*}
\<|W|^2x,;x\>
&\kern-6pt=\kern-6pt&
\|\kern1pt|W|x\|^2=\|\kern1ptU^*WTW^{-1}U|W|x\|^2                     \\
&\kern-6pt=\kern-6pt&
\|\kern1pt|W|TW^{-1}Wx\|^2=\|\kern1pt|W|Tx\|^2=\<T^*|W|^2Tx\,;x\>
\end{eqnarray*}
for every ${x\in\X}$, and hence ${T^*\!A'\,T\kern-1pt=A'}$ with 
${O\prec A'=|W|^2}$.
\end{proof}
\vskip-2pt

\vskip6pt
Let ${T\kern-1pt\in\kern-1pt\BX}$ and ${A\kern-1pt\in\kern-1pt\BX}$ be
arbitrary operators on a Hilbert space $\X.$ In ac\-cordance with
Proposition 4.1(a), $T$ was called an $A$-{\it isometry}\/ in \cite{Suc1} if
${T^*\!A\,T\kern-1pt=\kern-1ptA}$ for some ${A\ge O}.$ Similarly, $T$ was
called an $A$-{\it contraction}\/ in \cite{Suc1} if ${T^*\!A\,T\le A}$ for
some ${A\ge O}$ (see also \cite{Suc2, Suc3})$.$ Actually, $T$ is similar to a
contraction if and only if ${T^*\!A\,T\kern-1pt\le\kern-1pt A}$ some
${A\succ O}$ (see e.g., \cite[Corollary 1.8]{MDOT})$.$ Still along these
lines, if an operator $A$ (not necessarily nonnegative) is such that
${T^*\!A\,T=r(T)^2A}$ for some ${T\in\BX}$, where $r(T)$ stands for the
spectral radius of $T$, then $A$ was called $T$-{\it Toeplitz}\/ in
\cite{Ker5} (recall$:$ if $T$ is power bounded and ${\|T^n\|\not\to0}$, then
${r(T)=1}).$ For further applications of the new semi-inner product space
${(\X,\<\,\cdot\,;\cdot\,\>_A)}$ along different lines from those discussed
here see, e.g., \cite{ACG} and the references therein.

\vskip6pt
Similarity to an isometry is equivalent to the equation ${T^*\!A\,T=A}$ for
some ${A\succ O}$ (Proposition 4.1 (c,d))$.$ This is still equivalent to
some forms of power boundedness and power boundedness below (including
Ces\`aro means forms)$.$ These are brought together in the next proposition.

\vskip3pt\noi
\begin{proposition}
Let\/ $T$ be an operator on a Hilbert space\/ $\X.$ The following assertions
are pairwise equivalent.

\vskip4pt\noi
\begin{description}
\item{$\kern-9pt$\rm(a)$\kern1pt$}
$T$ is similar to an isometry.

\vskip6pt
\item{$\kern-9.5pt$\rm(b)$\kern1pt$}
$T$ is power bounded and power bounded below$:$ there exist\/
${\alpha,\beta\kern-1pt>\kern-1pt0}$ such that
$$
\alpha\|x\|\le\|T^kx\|\le\beta\|x\|
\;\;\hbox{for all}\;
k\ge0
\;\hbox{and every}\;
x\in\X.
$$

\vskip2pt
\item{$\kern-9pt$\rm(c)$\kern.5pt$}
There exist\/ ${\alpha,\beta>0}$ and\/ an invertible\/ ${R\in\BX}$ for which
$$
\alpha\|x\|^2\le\smallfrac{1}{n}{\sum}_{k=0}^{n-1}\|R\,T^kx\|^2\le\beta\|x\|^2
\;\;\hbox{for all}\;
n\ge1
\;\hbox{and every}\;
x\in\X.
$$

\vskip2pt
\item{$\kern-9pt$\rm(d)$\kern.5pt$}
There exist\/ ${\alpha,\beta>0}$ for which
$$
\alpha\|x\|^2\le\smallfrac{1}{n}{\sum}_{k=0}^{n-1}\|T^kx\|^2\le\beta\|x\|^2
\;\;\hbox{for all}\;
n\ge1
\;\hbox{and every}\;
x\in\X.
$$

\vskip2pt
\item{$\kern-9pt$\rm(e)$\kern.5pt$}
There exist\/ ${\alpha,\beta>0}$ and\/ an invertible\/ ${R\in\BX}$ for which
$$
\alpha\|x\|\le\|R\,T^kx\|\le\beta\|x\|
\;\;\hbox{for all}\;
k\ge0
\;\hbox{and every}\;
x\in\X.
$$
\end{description}
\end{proposition}

\vskip0pt\noi
\begin{proof}
Suppose (c) holds$.$ Since
$\frac{1}{n+1}\|R\,T^nx\|^2
\kern-1pt\le\kern-1pt\frac{1}{n+1}\kern-1pt\sum_{k=0}^n\|R\,T^kx\|^2
\kern-1pt\le\kern-1pt\beta\|x\|^2\kern-1pt$,
then
\vskip2pt\noi
\begin{eqnarray*}
\smallfrac{1}{n}\|R\,T^nx\|^2{\sum}_{k=0}^n\smallfrac{1}{k+1}
&\kern-6pt=\kern-6pt&
\smallfrac{1}{n}{\sum}_{k=0}^n\smallfrac{1}{k+1}\|R\,T^{k}T^{n-k}x\|^2
\le\smallfrac{1}{n}{\sum}_{k=0}^n\beta\|T^{n-k}x\|^2                        \\
&\kern-6pt=\kern-6pt&
\smallfrac{n+1}{n}\beta\smallfrac{1}{n+1}{\sum}_{k=0}^n\|R^{-1}R\,T^{k}x\|^2 
\le2\,\beta^2\|R^{-1}\|^2\|x\|^2,
\end{eqnarray*}
\vskip2pt\noi
and therefore
${\sup_n\big\|\big(\frac{1}{n}
\sum_{k=0}^n\frac{1}{k+1}\big)^{_{\scriptstyle\!\frac{1}{2}}}
R\,T^nx\big\|^{_{\scriptstyle 2}}\!<\kern-1pt\infty}$
for every ${x\in\X}.$ By the Banach--Steinhaus Theorem
$\,{\sup_n\frac{1}{n}\sum_{k=0}^n\frac{1}{k+1}\big\|R\,T^n\big\|^2
\!<\kern-1pt\infty}$
which implies ${\frac{1}{n}\|R\,T^n\|^2\!\to\kern-1pt0}$ (as
${\sum_{k=0}^n\frac{1}{1+k}\!\to\kern-1pt\infty}).$ Thus since
${\|T^{*n}R^*R\,T^n\|}={\|R\,T^n\|^2}$,
$$
\smallfrac{1}{n}\|T^{*n}R^*R\,T^n\|\to0.
$$
Now for each ${n\ge1}$ consider the Ces\`aro mean
$$
Q_n=\smallfrac{1}{n}{\sum}_{k=0}^{n-1}T^{*k}R^*R\,T^k.
$$
Since
${\|{Q_n}^\frac{1}{2}x\|^2\kern-2pt
=\kern-1pt\<Q_nx\,;x\>\kern-1pt
=\kern-1pt\frac{1}{n}\kern-2pt
\sum_{k=0}^{n-1}\kern-1pt\|R\,T^kx\|^2}\!$,
then
${\alpha\|x\|^2\kern-2pt
\le\kern-1pt\|{Q_n}^\frac{1}{2}x\|^2\kern-2pt\le\kern-1pt\beta\|x\|^2}\!$,
and hence $\{{Q_n}^\frac{1}{2}\kern-1pt\}$ is a bounded sequence of strictly
positive operators (and so is $\{Q_n\}$).

\vskip6pt\noi
(i) First suppose the Hilbert space $\X\kern-1pt$ is separable$.$ In this
case the bounded sequence $\{Q_n\}$ has a weakly convergent subsequence (see,
e.g., \cite[Theorem 5.70]{EOT})$.$ But the cone of nonnegative operators is
weakly closed in $\BX$ and so the weak limit of any weakly convergent
subsequence of $\{Q_n\}$ is again a nonnegative operator in $\BX.$ Let the
nonnegative ${Q\kern-1pt\in\kern-1pt\BX}$ be the weak limit of a weakly
convergent subsequence of $\{Q_n\}.$ Actually, $Q$ is strictly positive
because $\{Q_n\}$ is bounded below$.$ Since
$$
T^*Q_nT=Q_n+\smallfrac{1}{n}(R^*R-T^{*n}R^*R\,T^n)
$$
for each ${n\ge1}$, and since ${\smallfrac{1}{n}\|T^{*n}R^*R\,T^n\|\to0}$, we
get (again) the equation
$$
T^*Q\,T=Q.
$$
So $T$ is similar to an isometry (i.e., (a) holds) by Proposition 4.1(c). In
other words, with ${O\prec A=Q}$ the equation ${T^*\!A\,T=A}$ implies $T$
is similar to an isometry.

\vskip6pt\noi
(ii) Next suppose the Hilbert space $\X\kern-1pt$ is not separable$.$ Since
(c) holds, $T$ is nonzero$.$ Take an arbitrary
${0\kern-1pt\ne\kern-1ptx\kern-1pt\in\kern-1pt\X}$ and set
$\M_x\!=\kern-1pt{\span\big(\bigcup_{m,n\ge0}
\{T^mT^{*n}x\}\kern-1pt\cup\kern-1pt\{T^nT^{*m}x\}\big)^-}\!$
which is a separable nontrivial (closed) subspace of the (nonseparable)
Hilbert space $\X$ including $x$ and reducing $T.$ Then both
$T|_{\M_x}\kern-1pt$ and $(T|_{\M_x})^*\!=T^*|_{\M_x}\kern-1pt$ act on the
separable Hilbert space $\M_x.$ Thus since $T|_{\M_x}\kern-1pt$ satisfies (c)
--- because $T$ does --- then according to (i) $T|_{\M_x}\kern-1pt$ is
similar to an isometry$.$ Consider the collection
$\Im=$ $\{{\bigoplus\M_x}\!:{x\kern-1pt\in\kern-1pt\X}\}$ of all orthogonal
direct sums of these subspaces, which is partially ordered (in the inclusion
ordering) and is not empty (if
${0\kern-1pt\ne\kern-1pty\kern-1pt\in\kern-1pt
{\M_x}^{_{\scriptstyle\perp}}}\!$,
then ${\M_y\kern-1pt\sse\kern-1pt{\M_x}^{_{\scriptstyle\perp}}}\!$ because
${\M_x}^{_{\scriptstyle\perp}}$ \hbox{reduces} $T$ and so
${\M_x\kern-1pt\perp\kern-1pt\M_y}).$ Moreover, every chain in $\Im$ has an
upper bound in $\Im$ (the union of all orthogonal direct sums in a chain of
orthogonal direct sums in $\Im$ is again an orthogonal direct sum in $\Im$)$.$
Thus Zorn's Lemma ensures that $\Im$ has a maximal element, say
$\M={\bigoplus\M_x}$, which coincides with $\X$ (otherwise it would not be
maximal since ${\M\oplus\M^\perp}\!=\kern-1pt\X$)$.$ As $T|_{\M_x}x=Tx$, then
$T\kern-1pt={\bigoplus T|_{\M_x}}\kern-1pt$ on $\X\kern-1pt={\bigoplus\M_x}$
is similar to an isometry since each $T|_{\M_x}$ on $\M_x$ is similar to an
isometry according to item (i) above$.$ Thus again (a) holds$.$ Hence
$$
\rm (c)\;\limply\;(a).
$$
Now if (a) holds, then (as in the proof of (d) in Proposition 4.1) there
exists an invertible transformation ${W\in\B[\X,\Y]}$ for which
$\|WT^kW^{-1}y\|=\|y\|$ for every ${y\in\Y}$, equivalently,
$\|WT^kx\|=\|Wx\|$ for every ${x\in\X}$, for all ${k\ge0}.$ Therefore
$\|T^kx\|
\kern-1pt=\kern-1pt\|W^{-1}WT^kx\|
\kern-1pt\le\kern-1pt\|W^{-1}\|\kern1pt\|WT^kx\|
\kern-1pt=\kern-1pt\|W^{-1}\|\kern1pt\|Wx\|
\kern-1pt\le\kern-1pt\|W^{-1}\|\kern1pt\|W\|\kern1pt\|x\|$
and
$\|x\|
\kern-1pt=\kern-1pt\|W^{-1}Wx\|
\kern-1pt\le\kern-1pt\|W^{-1}\|\kern1pt\|Wx\|
\kern-1pt=\kern-1pt\|W^{-1}\|\kern1pt\|WT^kx\|
\kern-1pt\le\kern-1pt\|W^{-1}\|\kern1pt\|W\|\kern1pt\|T^kx\|$
so that
$\|W^{-1}\|^{-1}\|W\|^{-1}\|x\|
\kern-1pt\le\kern-1pt\|T^kx\|
\kern-1pt\le\kern-1pt\|W^{-1}\|\kern1pt\|W\|\kern1pt\|x\|$
for all $k$ for every $x.$ Hence
$$
\rm (a)\;\limply\;(b).
$$
This concludes the proof since
(b) $\limply$ (d) $\limply$ (c) $\Longleftarrow$ (e) $\Longleftarrow$ (b)
holds trivially.
\end{proof}

\vskip6pt
The positive numbers $\alpha$ and $\beta$ are constant with respect to $x$
but of course they may depend on $T.$ Parts of Proposition 4.2 have appeared
in \cite[Theorem 2]{KR}, \cite[Proposition 6.3]{RR}, \cite[Theorem 2]{Kan},
\cite[Proposition 1.15]{MDOT}, \cite[Theorem 1]{BS}, and
\cite[Corollary 4.2]{KD} (and parts of it --- e.g., (a)$\iff$(b) --- may
survive in a normed-space setting)$.$ For further conditions on similarity
to isometries along these lines see, e.g., \cite[Proposition 2.6]{Cas},
\cite[Theorem 2.1]{SS}.

\section{Contractions and the Equation ${T^*\!A\,T=A}$}

The next proposition is a classical result on Hilbert-space
contractions dating back to the early 1950s (see, e.g., \cite[Chapter 3]{MDOT}
and the references therein)$.$ It is based on a well-known result which
says$:$ {\it every monotone bounded sequence of Hilbert-space self-adjoint
operators converges strongly}\/. If $T$ is a Hilbert-space contraction, then
$\{T^{*n}T^n\}$ is a bounded monotone sequence of self-adjoint operators (in
fact a nonincreasing sequence of nonnegative operators) and so it converges
strongly to a nonnegative contraction $A$:
$$
\lim\,T^{*n}T^n=A
\quad\hbox{strongly}
$$
(i.e., ${\lim_n\|(T^{*n}T^n\!-A)x\|=0}$ for every $x).$ Such a nonnegative
contraction $A$ is usually refereed to as the {\it asymptotic limit}\/ of the
contraction $T.$ So if $T$ is a contraction, then the strong limit ${A\ge O}$
of $\{T^{*n}T^n\}$ (i.e., the asymptotic limit of $T$) defines a new
semi-inner product ${\<\,\cdot\,;\cdot\,\>_A}$ on $\X$ which becomes an inner
product if ${A>O}$, and in this case $T$ acts as an isometry on
${(\X,\<\,\cdot\,;\cdot\,\>_A)}$ by Proposition 4.1(b).

\vskip3pt\noi
\begin{proposition}
For every contraction\/ $T$ on a Hilbert space\/ $\X$ there exists a unique
nonnegative operator\/ $A$ on\/ $\X$ for which

\vskip6pt\noi
\begin{description}
\item{$\kern-9pt$\rm(a)$\kern2pt$}
$\kern-20pt$\centerline{
$T^{*n}T^n\sconv A$,
}$\kern-20pt$
\vskip2pt\noi
and so
$$
\<x\,;y\>_A={\lim}_n\<T^{*n}T^nx\,;y\>={\lim}_n\<T^nx\,;T^ny\>=\<Ax\,;y\>
$$
for every\/ ${x,y\in\X}$ or, equivalently,
$$
\|T^nx\|\to\|A^\frac{1}{2}x\|
\;\;\,\hbox{for every}\;\;
x\in\X,
$$
where\/ $\|A^\frac{1}{2}x\|=\|x\|_A=\|T^jx\|_A\,$ for every\/
${x\kern-1pt\in\kern-1pt\X}$ and every\/ ${j\kern-1pt\ge\kern-1pt0}$.
\end{description}

\vskip6pt\noi
Moreover\/$:$

\vskip4pt\noi
\begin{description}
\item{$\kern-9pt$\rm(b)$\kern2pt$}
${O\le A\le I}$\qquad
{\rm (i.e.,} $A$ is a nonnegative contraction on\/ $\X)$.

\vskip6pt
\item{$\kern-9pt$\rm(c)$\kern2pt$}
${T^*\!A\,T\!=\!A}.$\qquad
Equivalently,
\vskip4pt\noi
$\,\|A^\frac{1}{2}T^nx\|=\|A^\frac{1}{2}x\|\,$
for every\/ ${x\kern-1pt\in\!\X}$ and every\/ ${n\kern-1pt\ge\kern-1pt1}$.

\vskip6pt
\item{$\kern-9.5pt$\rm(d)$\kern2pt$}
$A\ne O\;\limply\;\|A\|=\|T\|=1$.

\vskip6pt
\item{$\kern-9.5pt$\rm(e)$\kern2pt$}
${A\,T=O\iff TA=O\iff A=O}$.

\vskip6pt
\item{$\kern-9pt$\rm(f)$\kern2pt$}
$A\,T=TA\iff\!A=A^2$.

\vskip6pt
\item{$\kern-9.5pt$\rm(g)$\kern2pt$}
${(I-A)\kern1ptT^n\sconv O}$
\quad and so \quad
${(I-A^\frac{1}{2})\kern1ptT^n\sconv O}$.

\vskip6pt
\item{$\kern-9.5pt$\rm(h)$\kern2pt$}
${\|A\kern1ptT^nx\|\to\|A^\frac{1}{2}x\|}\;\;$ for every\/ ${x\in\X}$.

\vskip6pt
\item{$\kern-8.5pt$\rm(i)$\kern2pt$}
$\N(A)=\big\{x\in\X\!:\,T^nx\to0\big\}$
\quad
$($so\/ $\;{T^n\!\sconv O}$ ${\!\iff\!}$ ${\!A=O})$.

\vskip6pt
\item{$\kern-9pt$\rm(j)$\kern2pt$}
${\N(I\kern-1pt-\kern-1ptA)\kern-1pt=\kern-1.5pt
\big\{x\kern-1pt\in\kern-1pt\X\!:\|T^nx\|\kern-1pt=\kern-1pt\|x\|
\;\;\forall n\kern-1pt\ge\kern-1.5pt1\big\}}$
\quad
$($so\/ $T$ is an isometry ${\!\iff\!\!}$ ${A\kern-1pt=\kern-1ptI})$.
\end{description}

\vskip6pt\noi
Furthermore,
$$
\hbox{$A$ is invertible $\iff$ $T$ is similar to an isometry.}
$$
\end{proposition}

\vskip0pt\noi
\begin{proof}
See, e.g., \cite[Propositions 3.1, 3.2 and 3.8]{MDOT} ---
see also \cite{Dur}, \cite{PV}, \cite{Dug}, \cite{KVP}.
\end{proof}
\vskip-2pt

\vskip3pt\noi
\begin{remark}
According to Proposition 5.1(i) the strong limit $A$ of $\{T^{*n}T^n\}$ for
a Hilbert-space contraction $T$ is positive if and only if $T$ is a
contraction of class $C_{1{\textstyle\cdot}}\!:$
$$
\hbox{$T$ is a $C_{1{\textstyle\cdot}}$-contraction}
\;\iff\;
\N(A)=\0
\;\iff\;
A>O.                                                        \leqno{\rm(a)}
$$
Since ${A^2\kern-2.5pt=\!A\!>\!O}$ implies ${A\!=\!I}\kern-1pt$, the above
equivalence and Proposition 5.1(j) ensures
$$
\hbox{$T$ is a $C_{1{\textstyle\cdot}}$-contraction with $A=A^2$}
\;\iff\;
\hbox{$T$ is an isometry}.                                  \leqno{\rm(b)}
$$
\end{remark}
\vskip-2pt

\vskip6pt
For a collection of properties of asymptotic limits for Hilbert-space
contractions see, for instance, \cite[Section I.10]{NF}, \cite[Section 6]{Fil},
\cite[Chapter 3]{MDOT}, \cite{Geh1}, \cite{Kub1}$.$ For the new inner product
${\<\cdot\,;\cdot\>_A}$ generated by the asymptotic limit $A$ of a
$C_{1{\textstyle\cdot}}$-contraction $T$ (i.e., for a positive $A$ or, in
particular, for a strictly positive $A$ as in Propositions 4.1(c) and 4.2)
see, for instance, \cite{Ker1}, \cite[Remark 3.9]{MDOT} and \cite{Kub0}.

\section{Power Bounded Operators and the Equation ${T^*\!A\,T=A}$}

The existence of Banach limits was established by Banach himself
\cite[p.21]{Ban} as a consequence of the Hahn--Banach Theorem$.$ Let
$\ell_+^\infty$ denote the Banach space of all complex-valued bounded
sequences equipped with its usual sup-norm$.$ A Banach limit is any
{\it bounded linear}\/ functional ${\vphi\!:\ell_+^\infty\!\to\CC}$ (i.e.,
${\vphi\in{\ell^\infty_+}^*}$, where ${\ell^\infty_+}^*$ is the dual of
${\ell^\infty_+}$) assigning a complex number to each complex-valued bounded
sequence which satisfies the following properties$.$ Take
${\{\xi_n\}\in\ell_+^\infty}$.

\vskip4pt\noi
\begin{description}
\item{$\kern-5.5pt$(o)}
$\;\vphi$ is linear\quad
(i.e., additive and homogenous),

\vskip6pt
\item{$\kern-4.5pt$(i)}
$\;\vphi$ is real\quad
(i.e., ${\vphi(\{\xi_n\})\in\RR}$ whenever $\{\xi_n\}$ is real-valued),

\vskip6pt
\item{$\kern-6pt$(ii)}
$\;\vphi$ is positive\quad
(i.e., ${0\le\vphi(\{\xi_n\})}$ whenever ${0\le\xi_n}$ for every $n$),

\vskip6pt
\item{$\kern-8pt$(iii)}
$\vphi$ is order-preserving\quad
(i.e., $\vphi(\{\xi_n\})\le\vphi(\{\upsilon_n\})$ if ${\xi_n\le\upsilon_n}$
in $\RR$ for every $n$),

\vskip6pt
\item{$\kern-8pt$(vi)}
$\;\vphi$ is backward-shift-invariant\quad
(i.e., $\vphi(\{\xi_{n+1}\})=\vphi(\{\xi_n\})$),

\vskip6pt
\item{$\kern-7pt$(v)}
$\;\liminf_n\xi_n\le\vphi(\{\xi_n\})\le\limsup_n\xi_n$ for every
real-valued sequence $\{\xi_n\}$,

\vskip6pt
\item{$\kern-9pt$(vi)}
$\;\vphi$ assigns to a convergent sequence its limit\quad
(i.e., ${\xi_n\conv\xi}\limply{\vphi(\{\xi_n\})=\xi}$)
\vskip2pt\noi
\kern4pt(in particular, ${\vphi(\{1,1,1,\dots\})=1}$),

\vskip6pt
\item{$\kern-12pt$(vii)}
$\;{\|\vphi\|=1}$.
\end{description}

\vskip6pt\noi
For existence of Banach limits see, e.g., \cite[Section III.7]{Con} or
\cite[Problem 4.66]{EOT})$.$ Moreover, for every Banach limit $\vphi$ there
exist Banach limits $\vphi_+$ and $\vphi_-$ such that
$$
{\lim}_n{\inf}_j\smallfrac{1}{n}{\sum}_{k=0}^{n-1}\xi_{k+j}
=\vphi_-(\{\xi_n\})
\le\vphi(\{\xi_n\})
\le\vphi_+(\{\xi_n\})
={\lim}_n{\sup}_j\smallfrac{1}{n}{\sum}_{k=0}^{n-1}\xi_{k+j}
$$
for an arbitrary real-valued sequence ${\{\xi_n\}\in\ell^\infty_+}$, where
$\vphi_-(\{\xi_n\})$ and $\vphi_+(\{\xi_n\})$ are the minimum and maximum
values of Banach limits at $\{\xi_n\}$, respectively
\cite[Theorem ($\beta,\gamma$)]{Such}; and for every
${\xi\in[\vphi_-(\{\xi_n\}),\vphi_+(\{\xi_n\})]}$ there exists a Banach limit
$\vphi'$ for which ${\vphi'(\{\xi_n\})=\xi}$ (see also \cite[(1.1)]{SeSu})$.$
Actually, all Banach limits coincide on a real-valued sequence if and only if
their shifted Ces\`aro means converge uniformly in the shifted parameter$.$
In other words, if $\{\xi_n\}$ is a real-valued sequence, then
$$
\vphi(\{\xi_n\})=\xi \;\;\hbox{for all Banach limits}\;\;
\vphi
\iff
{\lim}_n\smallfrac{1}{n}{\sum}_{k=0}^{n-1}\xi_{k+j}=\xi
\;\;\hbox{uniformly in $j$}.
$$
\cite[\!Theorem 1]{Lor} (see also \cite[Theorem ($\delta$)]{Such})$.$ We will
refer to the above displayed results as Lorentz characterizations$.$ Also, and
consequently, since $\vphi_-$ and $\vphi_+$ are Banach limits, then for every
real-valued sequence $\{\xi_n\}$,
$$
{\liminf}_n\xi_n
\le{\lim}_n{\inf}_j\smallfrac{1}{n}{\sum}_{k=0}^{n-1}\xi_{k+j}
\le{\lim}_n{\sup}_j\smallfrac{1}{n}{\sum}_{k=0}^{n-1}\xi_{k+j}
\le{\limsup}_n\xi_n.
$$

\vskip6pt
The Banach limit technique for power bounded operators discussed here is
well-known and has been applied quite often (see, e.g.,
\cite{Ker2, Ker3, Ker4} for applications along the lines considered here).

\vskip6pt
Suppose ${\!T\in\kern-1pt\BX}$ is a power bounded operator (i.e.,
${\sup_n\|T^n\|<\infty})$ acting on a Hilbert space
${(\X,\<\,\cdot\,;\cdot\,\>)}.$ Let ${\|\cdot\|\!:\X\!\to\kern-1pt\RR}$ be the
norm induced by the inner prod\-uct
${\<\,\cdot\,;\cdot\,\>\!:\X\!\times\!\X\!\to\CC}.$ Let
${\vphi\!:\ell_+^\infty\!\to\kern-1pt\CC}$ be an arbitrary Banach limit$.$
Since $T$ is power bounded, ${\{\<T^nx\,;T^ny\>\}\in\ell_+^\infty}$ for each
${x,y}$ in $\X.$ Thus set
$$
\<x\,;y\>_\vphi\!=\vphi(\{\<T^nx\,;T^ny\>\})=\vphi(\{\<T^{*n}T^nx\,;y\>\})
$$
for ${x,y}$ in $\X.$ Since every Banach limit is linear (which implies
$\vphi(\overline{\{\xi_n\}})=\overline\vphi(\{\xi_n\})$) and positive (i.e.,
${0\le\vphi(\{\xi_n\})}$ whenever ${0\le\xi_n}$ for every $n$), and since $T$
is linear, then it is readily verified that
${\<\cdot\,;\cdot\>_\vphi\!:\!\X\!\times\!\!\X\!\to\!\CC}$ is a semi-inner
product on $\X.$ Hence
$$
\|x\|_\vphi^2=\vphi(\{\|T^nx\|^2\})
$$
for every ${x\in\X}$ defines the seminorm ${\|\cdot\|_\vphi\!:\X\!\to\RR}$
induced by the semi-inner product ${\<\cdot\,;\cdot\>_\vphi}.$ (Even in this
case of a sequence of norms of powers of a power bounded operator, the squares
in the above identity cannot be omitted due to the nonmultiplicativity of
Banach limits)$.$ Since a Banach limit is order-preserving for real-valued
bounded sequences (i.e., $\vphi(\{\xi_n\})\le\vphi(\{\upsilon_n\})$ if
${\xi_n\le\upsilon_n}$ in $\RR$ for every $n$), and since
$\vphi(\{1,1,1,\dots\})=1$, then as $T$ is power bounded,
$$
\|x\|_\vphi\le{\sup}_n\|T^n\|\|x\|
$$
for every ${x\in\X}.$ Since Banach limits are backward-shift-invariant,
$$
\|Tx\|^2_\vphi=\vphi(\{\|(T^{n+1}x)\|^2\})
=\vphi(\{\|(T^nx)\|^2\})=\|x\|^2_\vphi.
$$
for every ${x\in\X}.$ The above setup leads to a generalization of
Proposition 5.1 from contractions to power bounded operators as in the
forthcoming Theorem 6.1.

\vskip3pt\noi
\begin{remark}
An important particular case$.$ If a power bounded operator $T$ is of class
$C_{1{\textstyle\cdot}}$ (i.e., ${T^nx\kern-1pt\notconv\kern-1pt0}$ for every
${0\kern-1pt\ne\kern-1ptx\kern-1pt\in\kern-1pt\X}$), then
${0\kern-1pt<\kern-1pt\liminf_n\|T^nx\|}$ for every
${x\kern-1pt\ne\kern-1pt0}$ (see, e.g., proof of Theorem 6.1(i) below)$.$
Any Banach limit $\vphi$ is such that
$\liminf_n\xi_n\kern-1pt\le\kern-1pt\vphi(\{\xi_n\})\le\limsup_n\xi_n$ for
a real-valued bounded sequence $\{\xi_n\}.$ Then
${0\kern-1pt<\kern-1pt\vphi(\{\|T^nx\|\})\kern-1pt=\kern-1pt\|x\|_\vphi}$
whenever ${x\kern-1pt\ne\kern-1pt0}.$ So the seminorm ${\|\cdot\|_\vphi}$
becomes a norm and consequently the semi-inner product
${\<\,\cdot\,;\cdot\,\>_\vphi}$ becomes an inner product$.$ Conversely, if
$T$ is not of class $C_{1{\textstyle\cdot}}$, then there is a nonzero
${x\kern-1pt\in\kern-1pt\X}$ for which ${T^nx\kern-1pt\conv\kern-1pt0}$ and
so ${\limsup_n\|T^nx\|\kern-1pt\to\kern-1pt0}$ which implies
${\vphi(\{\|T^nx\|^2|\})\kern-1pt=\kern-1pt0}$ so that
${\|x\|_\vphi\kern-1pt=\kern-1pt0}$, and the semi\-norm ${\|\cdot\|_\vphi}$
is not a norm$.$ Thus if ${\|\cdot\|_\vphi}\kern-1pt$ is a norm, then $T$ is
of class $C_{1{\textstyle\cdot}}.$ Hence
$$
\centerline{
${\<\,\cdot\,;\cdot\,\>_\vphi}$ is an inner product $\iff$ $T$ is a power
bounded of class $C_{1{\textstyle\cdot}}$.
}
$$
Thus if $T$ is a power bounded operator of class $C_{1{\textstyle\cdot}}$
(with respect to the original norm ${\|\cdot\|}$), then (since
${\|Tx\|_\vphi\kern-1pt=\kern-1pt\|x\|_\vphi}$ for every
${x\kern-1pt\in\kern-1pt\X}$ as we saw above) the norm ${\|\cdot\|_\vphi}$
makes $T$ into an isometry when acting on the inner product space
${(\X,\<\,\cdot\,;\cdot\,\>_\vphi)}$.
\end{remark}
\vskip-2pt

\vskip6pt
Proposition 5.1 can be extended from contractions to power bounded operators
(in particular, to power bounded operators of class $C_{1{\textstyle\cdot}}).$
Given a power bounded op\-erator $T$ and a Banach limit $\vphi$, there is a
unique nonnegative operator $A_\vphi$, referred to as the
$\vphi$-{\it asymptotic limit}\/ of $T\kern-1pt$, such that
${\<\cdot\,;\cdot\>_\vphi\kern-1pt=\kern-1pt\<A\cdot\,;\cdot\>}$ by
Lemma 3.1$.$ So
$$
\vphi(\{T^{*n}T^n\})=A_\vphi
\quad\hbox{weakly}
$$
(i.e., ${\vphi(\{\<T^{*n}T^nx\,;y\>\})=\<A_\vphi x\,;y\>}$ for every ${x,y}).$
The next theorem rounds up a collection of properties (either well-known or
not) of the $\vphi$-asymptotic limit $A_\vphi$ for a power bounded operator
$T$ into a unified statement$.$ Each assertion in Theorem 6.1 is written so as
to establish an injection from the items in Proposition 5.1 into homonymous
items in Theorem 6.1.

\vskip3pt\noi
\begin{theorem}
Let\/ ${T\kern-1pt\ne\kern-1ptO}$ be a power bounded operator on a Hilbert
space\/ ${(\X,\<\cdot\,;\cdot\>)}$ and let\/ ${\vphi\!:\ell_+^\infty\to\CC}$
be an arbitrary Banach limit\/$.$ Consider the semi-inner product\/
${\<\,\cdot\,;\cdot\,\>_\vphi}\!={\vphi(\{\<T^n\,\cdot\,;T^n\,\cdot\,\>\})}$
in\/ $\X$ generated by\/ $T$ and\/ $\vphi.$ Then there exists a unique
nonnegative operator\/ $A_\vphi$ on\/ $\X$ $($which depends on\/ $T$ and\/
$\vphi)$ such that\/
${\<\,\cdot\,;\cdot\,\>_\vphi}=\kern-1pt{\<A_\vphi\,\cdot\,;\,\cdot\>}$
and so

\vskip4pt\noi
\begin{description}
\item{$\kern-9pt$\rm(a)$\kern2pt$}
$\kern-20pt$\centerline{
$\<x\,;y\>_\vphi
=\,\vphi(\{\<T^nx\,;T^ny\>\})
=\,\vphi(\{\<T^{*n}T^nx\,;y\>\})
=\<A_\vphi x\,;y\>$
}$\kern-20pt$
\vskip6pt\noi
for every\/ ${x,y\in\X}$ or, equivalently,
$$
\vphi(\{\|T^nx\|^2\})=\|{A_\vphi}^{_{\!\scriptstyle\frac{1}{2}}}x\|^2
\;\;\,\hbox{for every}\;\;
x\kern-1pt\in\kern-1pt\X,
$$
where\/
$\|{A_\vphi}^{_{\!\scriptstyle\frac{1}{2}}}x\|=\|x\|_\vphi=\|T^jx\|_\vphi\,$
for every\/ ${x\kern-1pt\in\kern-1pt\X}$ and every\/
${j\kern-1pt\ge\kern-1pt0}$.
\end{description}

\vskip4pt\noi
Moreover\/$:$

\vskip3pt\noi
\begin{description}
\item{$\kern-9pt$\rm(b)$\kern2pt$}
$O\le A_\vphi\kern-.5pt\le\beta^2I\;$ with\/ $\;\beta={\sup}_n\|T^n\|\ne0$.
\vskip3pt\noi
Thus\/ $\|A_\vphi\|\le\beta^2$ and the identity\/ ${\|A_\vphi\|=\beta^2}$
may hold\/.

\vskip6pt
\item{$\kern-9pt$\rm(c)$\kern2pt$}
${T^*\kern-1ptA_\vphi\kern1ptT=A_\vphi}.$
\quad
Equivalently,
\vskip3pt\noi
$\,\|{A_\vphi}^{_{\!\scriptstyle\frac{1}{2}}}T^nx\|
=\|{A_\vphi}^{_{\!\scriptstyle\frac{1}{2}}}x\|\;$
for every\/ ${x\kern-1pt\in\!\X}$ and every\/ ${n\kern-1pt\ge\kern-1pt0}$.

\vskip6pt
\item{$\kern-9.5pt$\rm(d)$\kern2pt$}
$A_\vphi\ne O\;\limply\;1\le\|A_\vphi\|\;$ and\/ $\;1\le\|T\|$.

\vskip6pt
\item{$\kern-9.5pt$\rm(e)$\kern2pt$}
${A_\vphi\kern1ptT=O\iff TA_\vphi\kern-.5pt=O\iff A_\vphi\kern-.5pt=O}$.

\vskip6pt
\item{$\kern-9pt$\rm(f)$\kern2pt$}
$A_\vphi\kern1ptT=TA_\vphi\,\limply\;A_\vphi={A_\vphi}$$^{{\!\scriptstyle2}}$.

\vskip6pt\noi
Conversely,\quad
if\/ $A_\vphi={A_\vphi}$$^{{\!\scriptstyle2}}$, then
\vskip6pt\noi
{\rm (f$_1$)}
$\;\|(A_\vphi T^n-T^nA_\vphi)x\|^2\le(\beta^2\kern-1pt-1)\|A_\vphi x\|^2$
for all\/ $n$ and every\/ ${x\kern-1pt\in\kern-1pt\X}$,
\vskip3pt\noi
in particular,
$\;\|(A_\vphi T-TA_\vphi)x\|^2
\le(\|T\|^2\kern-1pt-1)\|A_\vphi x\|^2$
for every\/ ${x\kern-1pt\in\kern-1pt\X}$,
\vskip6pt\noi
{\rm (f$_2$)}
$\,\vphi(\{\|(A_\vphi T^n-T^nA_\vphi)x\|^2\})=0\,$
for every\/ ${x\kern-1pt\in\kern-1pt\X}$,
\vskip6pt\noi
{\rm (f$_3$)}
$\;{\|A_\vphi\|=1}\,$ whenever\/ $\kern1pt{O\ne A_\vphi}$.

\vskip6pt
\item{$\kern-9pt$\rm(g)$\kern2pt$}
$0\kern-1pt\le\kern-1pt\vphi(\{\|(I\kern-1pt-\kern-1ptA_\vphi)T^nx\|^2\})
\kern-1pt\le\kern-1pt
(\|A_\vphi\|^2\kern-1pt-\kern-1pt1)\,
\|{A_\vphi}^{_{\!\scriptstyle\frac{1}{2}}}x\|^2$
\qquad and
\vskip3pt\noi
$0\!\le\!\vphi(\{\|(I\!-\!{A_\vphi}^{_{\!\scriptstyle\frac{1}{2}}})T^nx\|^2\})
\!\le\!\|(I\kern-1pt+\kern-1pt{A_\vphi}^{_{\!\scriptstyle\frac{1}{2}}})^{-1}
\|^2(\|A_\vphi\|^2\!-\!1)\|{A_\vphi}^{_{\!\scriptstyle\frac{1}{2}}}x\|^2$
for any\/ ${x\kern-1pt\in\kern-1pt\X}$,
\vskip3pt\noi
which are both null if\/ ${\|A_\vphi\|=1}$ $($in particular if\/
$O\ne A_\vphi\kern-.5pt={A_\vphi}$$^{{\!\scriptstyle2}})\kern.5pt$
or\/ ${A_\vphi\kern-.5pt=O}$.

\vskip6pt
\item{$\kern-9.5pt$\rm(h)$\kern2pt$}
$\vphi(\{\|A_\vphi T^nx\|^2\})
=\|{A_\vphi}^{_{\!\scriptstyle\frac{1}{2}}}x\|^2
+\vphi(\{\|(I\!-\!A_\vphi)T^nx\|^2\})\,$
for every\/ ${x\kern-1pt\in\kern-1pt\X}$.

\vskip6pt
\item{$\kern-8.5pt$\rm(i)$\kern2pt$}
$\N(A_\vphi)
=\big\{{x\in\X\!:\,T^nx\to0}\big\}
=\big\{{x\in\X\!:\,\vphi(\{\|T^nx\|^2\})=0}\big\}$.
\vskip3pt\noi
Hence\/ $\,{\vphi(\{\|T^nx\|^2\})=0}$ for every ${x\in\X}\!$
${\!\iff\!}$ ${T^n\sconv O}$ $\!{\iff\!}$ ${A_\vphi\kern-1pt=O}$.

\vskip6pt
\item{$\kern-9pt$\rm(j)$\kern2pt$}
$\big\{x\in\X\!:\lim_n\|T^nx\|=\beta\|x\|\big\}$
\vskip2pt\noi
${\kern60pt}\sse
\N(\beta^2I-A_\vphi)
=\big\{x\in\X\!:\vphi(\{\|T^nx\|^2\})=\beta^2\|x\|^2\big\}$
\vskip2pt\noi
${\kern60pt}\sse
\big\{x\in\X\!:\|x\|\le{\liminf}_n\|T^nx\|
\le{\limsup}_n\|T^nx\|=\beta\|x\|\kern.5pt\big\}$.
\vskip3pt\noi
Hence\/ $\,{\vphi(\{\|T^nx\|^2\})=\beta^2\|x\|^2}\!$ for every\/
${x\kern-1pt\in\kern-1pt\X}$ ${\!\iff\!}$ ${\lim_n\|T^nx\|=\beta\|x\|}$ for
every\/ ${x\kern-1pt\in\kern-1pt\X}\!$ ${\!\iff\!}$
${A_\vphi\kern-.5pt=\beta^2I}$ ${\!\iff\!}$ ${A_\vphi\kern-.5pt=I}$
${\!\iff\!}$ $T$ is an isometry on\/ ${(\X,\<\,\cdot\,;\cdot\,\>)}$.
\end{description}

\vskip0pt\noi
Also,
\vskip0pt\noi
$$
\centerline{
${\<\,\cdot\,;\cdot\,\>_\vphi}$ is an inner product $\iff$ $T$ is of class
$C_{1{\textstyle\cdot}}\!\iff A_\vphi\kern-.5pt$\/ is positive.}
$$
\vskip3pt\noi
In this case\/ {\rm (i.e.,} if\/ ${A_\vphi\kern-.5pt>O})$,

\vskip3pt\noi
\begin{description}
\item{$\kern-9.5pt$\rm(k)$\kern2pt$}
$A_\vphi T\kern-1pt=\kern-1ptTA_\vphi
\!\iff\!A_\vphi\kern-1pt=\kern-1pt{A_\vphi}$$^{{\!\scriptstyle2}}\kern-1pt
\!\iff\!A_\vphi\kern-1pt=\kern-1ptI
\!\iff\!$
$\kern-1ptT$ is an isometry on\/ ${(\X,\<\,\cdot\,;\cdot\,\>)}$.
\end{description}

\vskip4pt\noi
Furthermore, the following assertions are pairwise equivalent.

\vskip2pt\noi
\begin{description}
\item{$\kern-0pt(1)\kern1pt$}
$A_\vphi$\/ is invertible\quad
{\rm(i.e.,} ${A_\vphi\kern-.5pt\succ O})$.

\vskip4pt
\item{$\kern-0pt(2)\kern0pt$}
The norms\/ ${\|\cdot\|_\vphi}$ and\/ ${\|\cdot\|}$ on\/ $\X$ are equivalent.

\vskip4pt
\item{$\kern-0pt(3)\kern1pt$}
$T$ on\/ ${(\X,\<\,\cdot\,;\cdot\,\>)}$ is similar to an isometry.

\vskip4pt
\item{$\kern-0pt(4)\kern1pt$}
$T$ on\/ ${(\X,\<\,\cdot\,;\cdot\,\>)}$ is power bounded below.
\end{description}
\end{theorem}

\vskip2pt\noi
{\it Proof}\/.
Let $T$ be a power bounded operator on a Hilbert space
${(\X,\<\,\cdot\,;\cdot\,\>)}.$ Thus the sequence $\{\<T^nx\,;T^ny\>\}$ is
bounded for every ${x,y\in\X}.$ Then consider the semi-inner product
${\<\,\cdot\,;\cdot\,\>_\vphi}=\vphi(\{\<T^n\,\cdot\,;T^n\,\cdot\>\})$ in
$\X$ generated by $T$ and a Banach limit $\vphi$.

\vskip4pt\noi
(a)
As ${(\X,\<\,\cdot\,;\cdot\,\>)}$ is a Hilbert space, an application of
Lemma 3.1 ensures the exist\-ence of a unique nonnegative operator $A_\vphi$
on $\X$ for each power bounded operator $T$ and each Banach limit $\vphi$,
the $\vphi$-asymptotic limit of $T$, such that
$$
\vphi(\{\<T^{*n}T^nx\,\;y\>\})
=\<x\,;y\>_\vphi=\<x\,;y\>_{A_\vphi}\!
=\<A_\vphi x\,;y\>
\;\;\,\hbox{for every}\;\;
x,y\in\X.
$$
The nontrivial part of the next equivalence follows by the polarization
identity$.$ The last identity is a consequence of the shift invariance
property for Banach limits:
$$
\|x\|_\vphi^2
=\|{A_\vphi}^{_{\!\scriptstyle\frac{1}{2}}}x\|^2\!
=\vphi(\{\|T^nx\|^2\})
=\vphi(\{\|T^nT^jx\|^2\})
=\|{A_\vphi}^{_{\!\scriptstyle\frac{1}{2}}}T^jx\|^2\!
=\|T^jx\|_\vphi^2
$$
for every ${x\kern-1pt\in\kern-1pt\X}$ and every ${j\kern-1pt\ge\kern-1pt0}$.

\vskip4pt\noi
(b) 
Now set $\beta\!=\!\sup_n\|T^n\|.$ Since
$\|{A_\vphi}^{_{\!\scriptstyle\frac{1}{2}}}x\|^2
\!=\!\vphi(\{\|T^nx\|^2\})
\!\le\!\sup_n\|T^nx\|^2
\!\le\!\beta^2\|x\|^2$
for every ${x\in\X}$ (because ${\|\vphi\|=1}$) we get
$$
\|A_\vphi\|=\|{A_\vphi}^{_{\!\scriptstyle\frac{1}{2}}}\|^2\le\beta^2.
$$
Also
${\<(A_\vphi\kern-1pt-\beta^2I)x\,;x\>}
\!=\!\|{A_\vphi}^{_{\!\scriptstyle\frac{1}{2}}}x\|^2\!-\!\beta^2\|x\|^2
\!\le\!(\|A_\vphi\|-\beta^2)\|x\|^2\le0$
for every ${x\kern-1pt\in\kern-1pt\X}$ by the above inequality$.$ Thus
the inequalities in (b) hold (since $A_\vphi$ is self-adjoint):
$$
O\le A_\vphi\kern-.5pt\le\beta^2I.
$$
(If $T={\rm shift}\{\beta,1,1,1,\dots\}$, then
$A_\vphi\kern-1pt={\rm diag}\{\beta^2\kern-1pt,1,1,1,\dots\}=T^{*n}T^n$
for all ${n\ge1}$.)

\vskip4pt\noi
(c)
By definition, ${\|x\|_{A_\vphi}\!=\|x\|_\vphi}$ and so
${\|Tx\|_{A_\vphi}\!=\|Tx\|_\vphi}.$ Since ${\|Tx\|_\vphi\!=\|x\|_\vphi}$
according to (a), then by (b) and Proposition 4.1(a)
$$
T^*\!A_\vphi\kern1ptT=A_\vphi.
$$
Equivalently, ${T^{*n}\!A_\vphi T^n=A_\vphi}$ for every ${n\ge1}$ by
induction, which means
$$
\|{A_\vphi}^{_{\!\scriptstyle\frac{1}{2}}}T^nx\|^2
=\<T^{*n}\!A_\vphi\kern1ptT^nx\,;x\>
=\<A_\vphi x\,;x\>
=\|{A_\vphi}^{_{\!\scriptstyle\frac{1}{2}}}x\|^2
$$
for every ${x\in\X}$ since $A_\vphi$ is nonnegative.

\vskip4pt\noi
(d)
Then
$\|{A_\vphi}^{_{\!\scriptstyle\frac{1}{2}}}x\|^2
\!=\kern-1pt\|{A_\vphi}^{_{\!\scriptstyle\frac{1}{2}}}T^nx\|^2
\!\le\|{A_\vphi}^{_{\!\scriptstyle\frac{1}{2}}}\|^2\|T^nx\|^2
\!=\kern-1pt\|A_\vphi\|\,\|T^nx\|^2.$
So (for a constant sequence)
$\|{A_\vphi}^{_{\!\scriptstyle\frac{1}{2}}}x\|^2
\!=\kern-1pt\vphi(\{\|{A_\vphi}^{_{\!\scriptstyle\frac{1}{2}}}x\|^2\})
\!\le\|A_\vphi\|\,\vphi(\{\|T^nx\|^2\})
\kern-1pt=\kern-1pt
\|A_\vphi\|\kern.5pt\|{A_\vphi}^{_{\!\scriptstyle\frac{1}{2}}}x\|^2$
for every ${x\in\X}$ by (a)$.$ Hence if there is ${x_0\kern-1pt\in\kern-1pt\X}$
for which
${{A_\vphi}^{_{\!\scriptstyle\frac{1}{2}}}x_0\kern-1pt\ne\kern-1pt0}$, then
${1\kern-1pt\le\kern-1pt\|A_\vphi\|}$:
$$
1\le\|A_\vphi\|
\quad\;\hbox{whenever}\;\quad
A_\vphi\kern-1pt\ne O.
$$
Since $A_\vphi=T^*\!A_\vphi\kern1ptT$, then
$\|A_\vphi\|\le\|A_\vphi\|\,\|T\|^2.$ Hence ${A_\vphi\ne O}$ implies
${1\le\|T\|}$.

\vskip4pt\noi
(e)
If ${A_\vphi=O}$, then ${A_\vphi\kern1ptT=TA_\vphi=O}$ trivially$.$ By (c),
${A_\vphi\kern1ptT=O}$ implies ${A_\vphi\kern-.5pt=O}.$ Finally, if\/
${TA_\vphi\kern-.5pt=O}$, then
$\vphi(\{\|T^nA_\vphi x\|\})\kern-1pt=\kern-1pt0$, and so
${\|{A_\vphi}^{_{\!\scriptstyle\frac{3}{2}}}x\|\kern-1pt=\kern-1pt0}$, for
every ${x\in\X}$, by (a)$.$ Thus ${A_\vphi\kern-.5pt=O}$ (by the Spectral
Theorem since by (b) $A_\vphi$ is nonnegative).

\vskip6pt\noi
(f)
Take ${x,y\kern-1pt\in\kern-1pt\X}.$ Since
${T^{*n}\!A_\vphi T^n}\!=\!A_\vphi$, then
${\vphi(\{\<T^{*n}\!A_\vphi T^nx\,;y\>\})}\kern-1pt=\kern-1pt
{\vphi(\{\<A_\vphi x\,;y\>\})}$
 $={\<A_\vphi x\,;y\>}$
(constant sequence)$.$ Also by (a)
${\vphi(\{\<T^{*n}T^nA_\vphi x\,;y\>\})}
\kern-1pt=\kern-1pt{\<{A_\vphi}^{_{\!\scriptstyle2}}x\,;y\>}.$
So
$$
A_\vphi T\kern-1pt=TA_\vphi\limply A_\vphi
\kern-1pt=\hbox{${A_\vphi}$$^{{\!\scriptstyle2}}$}.
$$
Conversely,
${\<A_\vphi T^nx\,;T^n\kern-1ptA_\vphi x\>}
={\|{A_\vphi}^{_{\!\scriptstyle\frac{1}{2}}}T^nx\|^2}
={\|{A_\vphi}^{_{\!\scriptstyle\frac{1}{2}}}x\|^2}$
by (c)$.$ If $A_\vphi\kern-.5pt={A_\vphi}$\!$^{{\!\scriptstyle2}}\kern-1pt$,
then ${{A_\vphi}^{_{\!\scriptstyle\frac{1}{2}}}\!=A_\vphi}$ by uniqueness of
the nonnegative square root, and hence
\begin{eqnarray*}
\|(A_\vphi T^n-T^n\kern-1ptA_\vphi) x\|^2
&\kern-6pt=\kern-6pt&
\|A_\vphi T^nx\|^2+\|T^n\kern-1ptA_\vphi x\|^2\!
-2\|{A_\vphi}^{_{\!\scriptstyle\frac{1}{2}}} x\|^2                      \\
&\kern-6pt=\kern-6pt&
\|T^n\kern-1ptA_\vphi x\|^2-\|A_\vphi x\|^2
\le(\beta^2-1)\|A_\vphi x\|^2
\end{eqnarray*}
for all $n$ and every ${x\in\X}.$ In particular, for ${n=1}$ we get for every
${x\in\X}$
$$
\|(A_\vphi T-TA_\vphi)x\|^2\le(\|T\|^2-1)\|A_\vphi x\|^2.
$$
Asymptotically, if
$A_\vphi\!=\kern-1pt{A_\vphi}$\!$^{{\!\scriptstyle2}}\kern-1pt$, then
${{A_\vphi}^{_{\!\scriptstyle\frac{1}{2}}}\!=\kern-1ptA_\vphi}.$ So we get
by (a) and the above identity
$$
\vphi(\{\|(A_\vphi T^n-T^n\kern-1ptA_\vphi)x\|^2\})
=\vphi(\{\|T^n\kern-1ptA_\vphi x\|^2\})-\|A_\vphi x\|^2
=\|{A_\vphi}^{_{\!\scriptstyle\frac{3}{2}}}x\|^2-\|A_\vphi x\|^2
=0
$$
for every ${x\in\X}.$ Moreover, if
$A_\vphi\kern-1pt={A_\vphi}$$^{{\!\scriptstyle2}}$, then $A_\vphi$ is an
orthogonal projection (since it is self-adjoint) and so ${\|A_\vphi\|=1}$
whenever ${A_\vphi\kern-1pt\ne O}$.

\vskip6pt\noi
(g)
Take an arbitrary ${x\in\X}.$
Again, by (c) we get
${\<T^nx\,;A_\vphi\,T^nx\>}
={\|{A_\vphi}^{_{\!\scriptstyle\frac{1}{2}}}T^nx\|^2}
={\|{A_\vphi}^{_{\!\scriptstyle\frac{1}{2}}}x\|^2}$
for all ${n\ge0}$, and
$\vphi(\{\|T^nx\|^2\})=\|{A_\vphi}^{_{\!\scriptstyle\frac{1}{2}}}x\|^2$
according to (a)$.$ Thus
\begin{eqnarray*}
0\le\vphi(\{\|(I-A_\vphi)T^nx\|^2\})
&\kern-6pt=\kern-6pt&
\vphi(\{\|T^nx\|^2\})+\vphi(\{\|A_\vphi T^nx\|^2\})
-2\vphi(\{\|{A_\vphi}^{_{\!\scriptstyle\frac{1}{2}}}x\|^2\})              \\
&\kern-6pt=\kern-6pt&
\vphi(\{\|A_\vphi T^nx\|^2\})-\|{A_\vphi}^{_{\!\scriptstyle\frac{1}{2}}}x\|^2
\le(\|A_\vphi\|^2-1)\|{A_\vphi}^{_{\!\scriptstyle\frac{1}{2}}}x\|^2.
\end{eqnarray*}
Since
${I-A_\vphi}
={(I+{A_\vphi}^{_{\!\scriptstyle\frac{1}{2}}})\,
(I-{A_\vphi}^{_{\!\scriptstyle\frac{1}{2}}})}$,
and since
${I+{A_\vphi}^{_{\!\scriptstyle\frac{1}{2}}}}$ is invertible with a bound\-ed
inverse (because ${A_\vphi\kern-.5pt\ge O}$), then
${I-{A_\vphi}^{_{\!\scriptstyle\frac{1}{2}}}}
={(I+{A_\vphi}^{_{\!\scriptstyle\frac{1}{2}}})^{-1}\kern1pt(I-A_\vphi)}$
and so
$$
0\le\vphi(\{\|(I-{A_\vphi}^{_{\!\scriptstyle\frac{1}{2}}})T^nx\|^2\})
\le\|(I+{A_\vphi}^{_{\!\scriptstyle\frac{1}{2}}})^{-1}\|^2
\vphi(\{\|(I-A_\vphi)T^nx\|^2\}).
$$

\vskip1pt\noi
(h)
This was proved above$:$
$\vphi(\{\|(I-A_\vphi)T^nx\|^2\})=
\vphi(\{\|A_\vphi T^nx\|^2\})-\|{A_\vphi}^{_{\!\scriptstyle\frac{1}{2}}}x\|^2$.

\vskip6pt\noi
(i)
Part of assertion (i) follows at once from (a) since
$\N(A_\vphi)=\N({A_\vphi}^{_{\!\scriptstyle\frac{1}{2}}}).$ Indeed,
$$
{T^nx\kern-.5pt\to\kern-.5pt0}
\limply
{\vphi(\{\|T^nx\|\})\kern-.5pt=\kern-.5pt0}
\!\!\iff\!\!
{\|{A_\vphi}^{_{\!\scriptstyle\frac{1}{2}}}x\|^2\!=\kern-.5pt0}
\!\!\iff\!\!
{x\kern-.5pt\in\kern-1pt\N({A_\vphi}^{_{\!\scriptstyle\frac{1}{2}}}\kern-1pt)}
\!\!\iff\!\!
{x\kern-.5pt\in\kern-1pt\N(A_\vphi)}.
$$
Conversely, suppose ${\beta\ge1}.$ (Otherwise ${T^nx\to0}$ for every
${x\kern-1pt\in\kern-1pt\X}$ since
${{\sup}_n\|T^n\|}={\beta\kern-1pt<\kern-1pt1}$ implies
${\|T^n\|\le\|T\|^n\le\beta^n\!\to0}.)$ If ${\vphi(\{\|T^nx\|\})=0}$ for some
Banach limit $\vphi$, then ${{\liminf}_n\|T^nx\|\kern-1pt=\kern-1pt0}$
(recall$:$ ${0\le{\liminf}_n\xi_n\!\le\vphi(\{\xi_n\})\le{\limsup}_n\xi_n}$
for ${\xi_n\kern-1pt\ge\kern-1pt0}).$ However, if ${\liminf}_n\|T^nx\|=0$,
then for every ${\veps\kern-1pt>\kern-1pt0}$ there is an integer $n_\veps$
such that ${\|T^{n_\veps}x\|<\veps}$, which implies
${\|T^nx\|\le\beta\|T^{n_\veps}x\|\kern-1pt<\kern-1pt\beta\kern.5pt\veps}$ for
all ${n\kern-1pt\ge\kern-1ptn_\veps}.$ Thus ${\|T^nx\|\to0}$.

\vskip6pt\noi
(j)
Take any ${0\kern-1pt\ne\kern-1ptx\kern-1pt\in\kern-1pt\X}$ (nonzero to
avoid trivialities)$.$ Since $\vphi$ is a Banach limit,
$$
\|T^nx\|\to\beta\|x\|
\;\limply\;
\vphi(\{\|T^nx\|^2\})=\beta^2\|x\|^2.
$$
According to (a),
$$
\!\vphi(\{\|T^nx\|^2\})=\beta^2\|x\|^2
\iff
\|{A_\vphi}^{_{\!\scriptstyle\frac{1}{2}}}x\|^2=\beta^2\|x\|^2
\iff
\<(\beta^2I-A_\vphi\kern-.5pt)x\,;x\>=0,
$$
and according to (b), since
${\N((\beta^2I-A_\vphi\kern-.5pt)^\frac{1}{2}\kern-1pt)
=\N(\beta^2I-A_\vphi\kern-.5pt)}$,
$$
\<(\beta^2I-A_\vphi\kern-.5pt)x\,;x\>=0
\iff
\|(\beta^2I-A_\vphi\kern-.5pt)^\frac{1}{2}x\|=0
\iff
x\in\N(\beta^2I-A_\vphi\kern-.5pt).
$$
Conversely, since $\vphi(\{\|T^nx\|^2\})=\beta^2\|x\|^2$ means
${\vphi(\{\beta^2\|x\|^2\kern-1pt-\|T^nx\|^2\})=0}$, and since
${0\le\beta^2\|x\|^2\kern-1pt-\|T^nx\|^2}$ for every
$n$ because ${{\sup}_n\|T^nx\|\le\beta\|x\|}$, then we get
$$
\vphi(\{\|T^nx\|^2\})=\beta^2\|x\|^2
\;\limply\;
{\liminf}_n\{\beta^2\|x\|^2-\|T^nx\|^2\}=0.
$$
However, recalling again that ${0\le\|T^nx\|\le\beta\|x\|}$ for every $n$,
$$
{\liminf}_n\{\beta^2\|x\|^2-\|T^nx\|^2\}=0
\iff
{\limsup}_n\|T^nx\|=\beta\|x\|.
$$
Moreover, for ${x\kern-1pt\ne\kern-1pt0}$ and since
${\beta\kern-1pt>\kern-1pt0}$ (as ${T\kern-1pt\ne\kern-1ptO}$),
$$
{\limsup}_n\|T^nx\|=\beta\|x\|
\;\limply\;
\beta\ge1
$$
(indeed, if ${\beta\kern-1pt<\kern-1pt1}$, then
$\beta={{\limsup}_n\kern-1pt\frac{\|T^nx\|}{\|x\|}}
\kern-1pt\le\kern-.5pt
{{\limsup}_n\kern-1pt\sup_{x\ne0}\kern-2pt\frac{\|T^nx\|}{\|x\|}}
\kern-1pt=\kern-.5pt
{{\limsup}_n\|T^n\|}
\kern-1pt\le\kern-.5pt
{{\limsup}_n\|T\|^n}
\kern-1.5pt\le\kern-.5pt
{{\limsup}_n\beta^n}
\kern-1.5pt=\kern-.5pt
{{\lim}_n\beta^n}
\kern-1.5pt=\kern-.5pt
0).$
Also, since ${\|T^{n+m}x\|\le\beta\|T^nx\|}$ for each ${m,n\ge0}$, then 
${\limsup}_n\|T^nx\|={\limsup}_n\|T^{n+m}x\|\le\beta\,{\liminf}_m\|T^mx\|$,
and so
$$
{\limsup}_n\|T^nx\|=\beta\|x\|
\;\limply\;
\|x\|\le{\liminf}_m\|T^mx\|.
$$
Finally, by the above implications and equivalences, if
$\lim_n\|T^nx\|=\beta\|x\|$ for every ${x\kern-1pt\in\kern-1pt\X}$, then
$\vphi(\{\|T^nx\|^2\})=\beta^2\|x\|^2$ for every ${x\kern-1pt\in\kern-1pt\X}$,
which means $A_\vphi\kern-1pt=\beta^2I.$ But this implies ${T^*T=I}$ by (c)
(since ${\beta\ne0}$ whenever ${T\ne O})$, which in turn implies
${A_\vphi\kern-1pt=I}$ by (a) (i.e.,
${\<A_\vphi x\,;y\>}=\vphi(\{\<T^{*n}T^nx\,;y\>\})$ for every
${x,y\kern-1pt\in\kern-1pt\X}).$ However, if ${A_\vphi\kern-1pt=I}$, then
${\|Tx\|=\|x\|}$ for every ${x\kern-1pt\in\kern-1pt\X}$ by (c), which means
$T$ is an isometry on ${(\X\kern-1pt,\kern-1pt\<\,\cdot\,;\cdot\,\>)}$, and
we are back to ${\lim_n\|T^nx\|=\beta\|x\|}$ for every
${x\kern-1pt\in\kern-1pt\X}$ with ${\beta=1}$.

\vskip6pt
As we saw in Remark 6.1, ${\<\,\cdot\,;\cdot\,\>_\vphi}$ is an inner product
if and only if $T$ is a power bounded of class $C_{1{\textstyle\cdot}}$, and
the semi-inner product
${\<\,\cdot\,;\cdot\,\>_\vphi}={\<A_\vphi\,\cdot\,;\cdot\,\>}$ is an inner
product (i.e., the seminorm
${\|\cdot\|_\vphi}={\|{A_\vphi}^{_{\!\scriptstyle\frac{1}{2}}}\cdot\|}$ is a
norm) if and only if $\N(A_\vphi)=\0$, which means the nonnegative $A_\vphi$
is positive$.$ Thus from now on suppose ${A_\vphi\!>O}$.

\vskip6pt\noi
(k)
If ${A_\vphi T\kern-1pt=TA_\vphi}$ then
$A_\vphi\!=\kern-1pt{A_\vphi}$$^{{\!\scriptstyle2}}$ by (f) and so
$A_\vphi\kern-1pt$ is an orthogonal projection (since it is a seif-adjoint
idempotent) which implies ${A_\vphi\kern-1pt=I}$ (because
${A_\vphi\!>O})$, and hence ${A_\vphi T\kern-1pt=TA_\vphi}$
trivially$.$ Therefore if ${A_\vphi\!>O}$,
$$
A_\vphi T\kern-1pt=TA_\vphi
\;\limply\; A_\vphi\kern-1pt=\hbox{${A_\vphi}$$^{{\!\scriptstyle2}}$}
\;\limply\; A_\vphi\kern-1pt=I
\;\limply\; A_\vphi T\kern-1pt=TA_\vphi.
$$
But ${A_\vphi\kern-1pt=I}$ if and only if $T$ is an isometry, as we saw in
the proof of item (j).

\vskip6pt
Since $T$ is power bounded, to prove assertions (1) to (4) proceed as
follows$.$
$$
\rm
(1)\iff(2),\quad
(3)\iff(4),\quad
and\quad
(1)\;\,\limply\;\,(3)
$$
by Proposition 3.1, Proposition 4.2(a,b), and Proposition 4.1(c),
respectively$.$ Conversely, as $T$ is power bounded, if (4) holds, then
$\alpha\|x\|\le\|T^nx\|\le\beta\|x\|$ for all ${n\ge0}$, and so
$\alpha\|x\|\le\vphi(\{\|T^nx\|\})\le\beta\|x\|$, for every ${x\in\X}.$
Since $\vphi(\{\|T^nx\|^2\})=\|x\|_\vphi^2$ by (a), then
${\alpha^2\|x\|^2\le\|x\|_\vphi^2\le\beta^2\|x\|^2}$ for every ${x\in\X}$
and so (2) holds$.$ Thus
$$
(4)\limply(2). \eqno{\qed}
$$
\vskip-2pt

\vskip3pt\noi
\begin{remark}
If $T$ is a contraction (equivalently, if ${\beta\kern-1pt\le\kern-1pt1}$),
then Theorem 6.1 is reduced to Proposition 5.1, and $\{T^{*n}T^n\}$ converges
strongly (thus weakly) to ${A_\vphi\kern-1pt=A}$ for every Banach limit
$\vphi$, with ${\|A\|=1}$ or ${\|A\|=0}.$ For a
$C_{1{\textstyle\cdot}}$-contraction, ${\|A\|=1}$.
\end{remark}

\vskip6pt
Such a combined procedure (of using Lemma 3.1 together with an inner product
generated by a power bounded operator and a Banach limit) seems to have been
originated in the celebrated Nagy's 1947 paper \cite{Nag} (see also
\cite[Section II.5]{NF})$.$ Subsequent applications of it appear, for
instance, in \cite{Ker1, Ker2, Ker3} and, recently, in
\cite{Geh2, Geh3, KD}$.$ Proposition 4.2(a,b), however, supplies an
elementary and straightforward proof of Nagy's result as follows.

\vskip3pt\noi
\begin{corollary}
\cite{Nag}
On a Hilbert space, an invertible power bounded operator with a power bounded
inverse is similar to a unitary operator --- the converse is trivial\/.
\end{corollary}

\vskip0pt\noi
\begin{proof}
Let ${T\in\BX}$ and ${T^{-1}\in\BX}$ be power bounded$.$ Thus there exist
real con\-stants ${0<\alpha\le1}$ and ${1\le\beta}$ for which
${\|T^n\|\le\beta}$ and ${\|T^{-n}\|\le\alpha^{-1}}$ for all ${n\ge0}.$ So
${\alpha\|x\|}\le{\|T^{-n}\|^{-1}\|x\|}\le{\|T^nx\|\le\beta\|x\|}$ for all $n$
and every $x.$ Hence $T$ is similar to an isometry by Proposition 4.2(a,b)$.$
Since $T$ is invertible, then so is the isometry similar to it: an invertible
Hilbert-space isometry means a unitary operator.
\end{proof}

\section{Ces\`aro Means and the Equation ${T^*\!A\,T=A}$}

A word on terminology$.$ An $\X$-valued sequence in an arbitrary normed
space $\X$ is called Ces\'aro convergent if its sequence of arithmetic
means (referred to as Ces\'aro means) converges in $\X$, whose limit is
called Ces\'aro limit.

\vskip6pt
Banach limits have been related to Ces\`aro means since the very beginning
\cite{Lor}, and Ces\`aro means are naturally linked to the Ergodic Theorem
for power bounded operators$.$ If a sequence $\{Q_n\}$ of Ces\`aro means
$Q_n=\smallfrac{1}{n}{\sum}_{k=0}^{n-1}T^{*k}T^k$ for an operator $T$
converges (either weakly, strongly, or uniformly), then its limit $Q$ (if it
exists) has been refereed to as the {\it Ces\`aro asymptotic limit}\/ of $T$
(see \cite{Geh2})$.$ As is well-known, the strong limit $Q$ always exists for
contractions and coincides with the asymptotic limit $A$$:$ for a contraction
$T$ the sequence of Ces\`aro means $\{Q_n\}$ converges strongly to ${Q=A}.$
An elementary quick proof is readily obtained as follows.

\vskip3pt\noi
\begin{proposition}
$\kern-2.5pt$For$\kern-.5pt$ a$\kern-.5pt$ Hilbert-space$\kern-.5pt$
contraction$\kern-.5pt$ the$\kern-.5pt$ sequence$\kern-.5pt$ of$\kern-.5pt$
Ces\`aro$\kern-.5pt$ mean$\kern-.5pt$s con\-verges$\kern-.5pt$
strongly$\kern-.5pt$ and$\kern-.5pt$ the$\kern-.5pt$ Ces\`aro$\kern-.5pt$
asymptotic$\kern-.5pt$ limit$\kern-.5pt$ coincides$\kern-.5pt$
with$\kern-.5pt$ the$\kern-.5pt$ asymptotic$\kern-.5pt$ limit\/$:$
$$
\|T\|\le1\quad\limply\quad Q_n=\smallfrac{1}{n}{\sum}_{k=0}^{n-1}T^{*k}T^k
\sconv Q=A.
$$
\end{proposition}

\vskip0pt\noi
\begin{proof}
If $T$ is a contraction, then the sequence $\{T^{*n}T^n\}$ converges strongly
to $A$ by Proposition 5.1$.$ Since ${x_n\to x}$ implies
${\frac{1}{n}\sum_{k=1}^n x_k\to x}$ for any normed-space-valued sequence
$\{x_n\}$, then ${T^{*n}T^n\sconv A}$ implies ${Q_n\sconv Q=A}$.
\end{proof}
\vskip-2pt

\vskip6pt
As before, the next theorem brings together scattered properties (again,
either well-known --- e.g., \cite[Theorems 2.5, 2.6 and Proposition 5.1]{Geh2}
--- or not) of Ces\`aro asymptotic limits $Q$ into a unified statement$.$
Some parts in the proof behave similarly to their equivalent in the proof of
Theorem 6.1, as expected; some other parts require an independent and
different approach$.$ Each assertion in Theorem 7.1 below is written so as to
establish a bijection with the items in Theorem 6.1 and so, by transitivity,
it establishes an injection from the items in Proposition 5.1 into homonymous
items in Theorems 6.1 and 7.1.

\vskip3pt\noi
\begin{theorem}
Let\/ ${O\ne T\kern-1pt\in\kern-1pt\BX}$ be a Hilbert-space operator$.$ For
each positive in\-teger\/ $n$ consider the Ces\`aro mean
$$
Q_n={\smallfrac{1}{n}{\sum}_{k=0}^{n-1}T^{*k}T^k}
$$
in $\BX.$ Suppose the sequence\/ $\{Q_n\}$ converges weakly to\/
${Q\kern-1pt\in\kern-1pt\BX}.$ That is, suppose

\vskip2pt\noi
\begin{description}
\item{$\kern-9pt$\rm(a)$\kern2pt$}
$Q_n\wconv Q.$
\quad
Equivalently,
\vskip4pt\noi
$\,\|{Q_n}^\frac{1}{2}x\|\to\|Q^\frac{1}{2}x\|\;$
for every\/ ${x\kern-1pt\in\!\X}$.
\end{description}

\vskip4pt\noi
Then in this case\/$:$

\vskip2pt\noi
\begin{description}
\item{$\kern-9pt$\rm(b)$\kern2pt$}
${O\le Q}\;$ and, if\/ ${\sup_n\|T^n\|=\beta}$ $($so that\/ ${\beta\ne0})$,
then\/ $\;{Q\le\beta^2}$.
\vskip4pt\noi
If\/ $T$ is power bounded, then\/ $\|Q\|\le\beta^2$ and the identity\/
${\|Q\|=\beta^2}\kern-1pt$ may hold\/.

\vskip4pt
\item{$\kern-9pt$\rm(c)$\kern2pt$}
${T^*Q\,T=Q}.$
\quad
Equivalently,
\vskip4pt\noi
$\,\|Q^{_{\!\scriptstyle\frac{1}{2}}}T^nx\|
=\|Q^{_{\!\scriptstyle\frac{1}{2}}}x\|\;$
for every\/ ${x\kern-1pt\in\!\X}$ and every\/ ${n\kern-1pt\ge\kern-1pt0}.$
\quad
Therefore
\vskip4pt
\item{$\kern5pt$}
$\|{Q_n}^\frac{1}{2}T^jx\|^2\!
=\smallfrac{1}{n}{\sum}_{k=0}^{n-1}\|T^{k+j}x\|^2\to\|Q^\frac{1}{2}x\|^2\,$
for every\/ ${x\in\X}$ and every\/ $j\ge0$.

\vskip4pt
\item{$\kern-9.5pt$\rm(d)$\kern2pt$}
$Q\ne O\;\limply\;1\le\|Q\|\;$ and\/ $\;1\le\|T\|$.

\vskip4pt
\item{$\kern-9.5pt$\rm(e)$\kern2pt$}
${Q\kern1ptT=O\iff TQ=O\iff Q=O}$.

\vskip4pt
\item{$\kern-9pt$\rm(f)$\kern2pt$}
${Q\kern1ptT=TQ\,\limply\;Q=Q^2}$.

\vskip6pt\noi
Conversely,\quad
if\/ ${Q=Q^2}$, then
\vskip4pt\noi
{\rm (f$_1$)}
$\;\|(Q\kern1ptT^n-T^nQ)x\|^2\le(\sup_n\|T^n\|^2\kern-1pt-1)\|Qx\|^2$
for all\/ $n$ and every\/ ${x\kern-1pt\in\kern-1pt\X}$,
\vskip4pt\noi
in particular,
$\;\|(Q\kern1ptT-TQ)x\|^2\le(\|T\|^2\kern-1pt-1)\|Qx\|^2$
for every\/ ${x\kern-1pt\in\kern-1pt\X}$,
\vskip4pt\noi
{\rm (f$_2$)}
$\,\|(Q\kern1ptT^n-T^nQ)x\|\to0\,$
for every\/ ${x\kern-1pt\in\kern-1pt\X}$,
\vskip4pt\noi
{\rm (f$_3$)}
$\;{\|Q\|=1}\,$ whenever\/ $\kern1pt{O\ne Q}$.

\vskip4pt
\item{$\kern-9pt$\rm(g)$\kern2pt$}
$\|(I-Q){Q_n}^\frac{1}{2}x\|^2\le(\|Q\|^2-1)\,\|{Q_n}^\frac{1}{2}x\|^2
\,\to\;(\|Q\|^2-1)\,\|Q^\frac{1}{2}x\|^2$
\quad
and
\vskip4pt\noi
$\,\|(I-Q^\frac{1}{2}){Q_n}^\frac{1}{2}x\|^2
\le\|(I+Q^\frac{1}{2})^{-1}\|^2(\|Q\|^2-1)\,\|{Q_n}^\frac{1}{2}x\|^2$
for every\/ ${x\kern-1pt\in\kern-1pt\X}$,
\vskip4pt\noi
which are both null if\/ ${\|Q\|=1}$ or\/ asymptotically null if\/ ${Q=O}$.

\vskip4pt
\item{$\kern-9.5pt$\rm(h)$\kern2pt$}
$\|Q\kern1pt{Q_n}^\frac{1}{2}x\|^2
=\|(2Q-I)^\frac{1}{2}{Q_n}^\frac{1}{2}x\|^2
+\|(I-Q)\kern1pt{Q_n}^\frac{1}{2}x\|^2$
for every\/ ${x\kern-1pt\in\kern-1pt\X}$,

\vskip4pt
\item{$\kern-8.5pt$\rm(i)$\kern2pt$}
If\/ $T$ is power bounded, then\/
\vskip2pt\noi
$\N(Q)
=\big\{{x\in\X\!:\,T^nx\to0}\big\}
=\big\{x\in\X\!:\,\frac{1}{n}{\sum}_{k=0}^{n-1}\|T^kx\|\to0\big\}$.
\vskip4pt\noi
Hence\/ $\;{T^n\!\sconv O}$ ${\!\iff\!}$ ${Q=O}$ and\/ $T$ is power bounded.

\vskip4pt
\item{$\kern-9pt$\rm(j)$\kern2pt$}
If\/ $\,\sup_n\|T^n\|=\beta$ $\;($so that\/ ${\beta\ne0})\;$ then\/
\vskip2pt\noi
$\big\{x\in\X\!:\lim_n\|T^nx\|=\beta\|x\|\kern.5pt\big\}$
\vskip2pt\noi
${\kern60pt}\sse
\N(\beta^2I-Q)
=\big\{x\in\X\!:\lim_n\smallfrac{1}{n}{\sum}_{k=0}^{n-1}\|T^kx\|^2
=\beta^2\|x\|\kern.5pt\big\}$
\vskip2pt\noi
${\kern60pt}\sse
\big\{x\in\X\!:\|x\|\le{\liminf}_n\|T^nx\|\le{\limsup}_n\|T^nx\|
=\beta\|x\|\kern.5pt\big\}$.
\vskip4pt\noi
Hence\/ ${\lim}_n\smallfrac{1}{n}{\sum}_{k=0}^{n-1}\|T^kx\|^2=\beta^2\|x\|$
for every\/ ${x\kern-1pt\in\kern-1pt\X}$ ${\!\iff\!}$
${\lim_n\|T^nx\|=\beta\|x\|}$ for every\/
${x\in\X}\!$ ${\!\iff\!}$ ${Q=\beta^2I}$ ${\!\iff\!}$
${Q=I}$ ${\!\!\iff\!\!}$ $T$ is an isometry on
${(\X,\<\,\cdot\,;\cdot\,\>)}$.
\end{description}

\vskip4pt\noi
Also if\/ $T$ is power bounded, then
$$
\centerline{
${\<\,\cdot\,;\cdot\,\>_Q}$ is an inner product $\iff$ $T$ is of class
$C_{1{\textstyle\cdot}}\!\iff Q$\/ is positive.}
$$
In the case of\/ ${Q>O}$ we get

\vskip2pt\noi
\begin{description}
\item{$\kern-9.5pt$\rm(k)$\kern2pt$}
$Q\kern1ptT=TQ
\!\iff\!Q=Q^2\kern-1pt
\!\iff\!Q=I
\!\iff\!$
$\kern-1ptT$ is an isometry on\/ ${(\X,\<\,\cdot\,;\cdot\,\>)}$.
\end{description}

\vskip4pt\noi
Furthermore, the following assertions are pairwise equivalent.

\vskip2pt\noi
\begin{description}
\item{$\kern-0pt(1)\kern1pt$}
$Q$\/ is invertible\quad
{\rm(i.e.,} ${Q\succ O})$.

\vskip4pt
\item{$\kern-0pt(2)\kern0pt$}
The norms\/ ${\|\cdot\|_Q}$ and\/ ${\|\cdot\|}$ on\/ $\X$ are equivalent.

\vskip4pt
\item{$\kern-0pt(3)\kern1pt$}
$T$ on\/ ${(\X,\<\,\cdot\,;\cdot\,\>)}$ is similar to an isometry.

\vskip4pt
\item{$\kern-0pt(4)\kern1pt$}
$T$ on\/ ${(\X,\<\,\cdot\,;\cdot\,\>)}$ is power bounded below.
\end{description}
\end{theorem}

\vskip0pt\noi
\begin{proof}
Let ${T\kern-1pt\in\kern-1pt\BX}$ be an operator acting on a Hilbert space
$\X$, take the sequence $\{T^{*n}T^n\}$ of nonnegative operators in $\BX$,
and consider the Ces\`aro mean
$$
Q_n={\smallfrac{1}{n}{\sum}_{k=0}^{n-1}T^{*k}T^k}
\quad\;\hbox{for every}\;\;n\ge1
$$
associated with $\{T^{*n}T^n\}.$ Suppose the $\BX$-valued sequence $\{Q_n\}$
of nonnegative operators converges weakly to the Ces\`aro asymptotic limit
${Q\kern-1pt\in\kern-1pt\BX}$ of $T$.

\vskip6pt\noi
(a)
Since the class of nonnegative operators is weakly closed in $\BX$, then the
weak limit $Q$ is nonnegative, and hence ${Q_n\wconv Q}$ is equivalent to
$$
\|{Q_n}^\frac{1}{2}x\|^2
=\<Q_nx\,\;x\>
\to
\<Qx\,\;x\>
=\|Q^\frac{1}{2}x\|^2
\quad\;\hbox{for every}\;\;x\in\X.
$$
\vskip-4pt

\vskip6pt\noi
(b)
By (a), ${O\le Q}$. If\/ $\,{\sup_n\|T^n\|=\beta}$, then\/
$Q_n\le\smallfrac{1}{n}I
+{\smallfrac{1}{n}{\sum}_{k=1}^{n-1}\beta^2I}\kern-1pt\to\beta^2I$
(and if\/ ${T\kern-1pt={\rm shift}\{\beta,1,1,1,\dots\}}$, then\/
${Q=T^{*n}T^n\kern-1pt={\rm diag}\{\beta^2\kern-1pt,1,1,1,\dots\}})\,$
for every ${n\kern-1pt\ge\kern-1pt1}$.

\vskip6pt\noi
(c)
Since $Q_n={\smallfrac{1}{n}{\sum}_{k=0}^{n-1}T^{*k}T^k}\!$, then
(compare with the proof of Proposition 4.2)
$$
T^*Q_nT=Q_{n+1}+\smallfrac{1}{n}(Q_{n+1}-I)
$$
for each ${n\kern-1pt\ge\kern-1pt1}.$ If
${\<Q_nx\,;x\>\kern-1pt\to\kern-1pt\<Qx\,;x\>}$ for every
${x\kern-1pt\in\kern-1pt\X}$, then $\kern-.5pt\{Q_n\}\kern-.5pt$ is bounded
and so
$$
T^*Q\,T=Q.
$$
By induction, ${T^{*n}Q\,T^n\kern-1pt=Q}$ for all
${n\kern-1pt\ge\kern-1pt0}.$ As ${O\le Q}$ and ${O\le Q_n}$ for
${n\kern-1pt\ge\kern-1pt1}$, then (i)
$$
\|Q^\frac{1}{2}T^nx\|^2
=\<Q\kern1ptT^nx\,;T^nx\>
=\<T^{*n}Q\kern1ptT^nx\,;x\>
=\<Qx\,;x\>
=\|Q^\frac{1}{2}x\|
$$
for every ${x\kern-1pt\in\kern-1pt\X}$ and every ${n\kern-1pt\ge\kern-1pt0}$,
and (ii) ${T^{*j}Q_nT^j\kern-1pt\wconv Q}$ for every
${j\kern-1pt\ge\kern-1pt0}$, and so
$$
\|{Q_n}^\frac{1}{2}T^jx\|^2
=\<Q_nT^jx\,;T^jx\>
=\smallfrac{1}{n}{\sum}_{k=0}^{n-1}\<T^{*k+j}T^{k+j}x\,;x\>
=\smallfrac{1}{n}{\sum}_{k=0}^{n-1}\|T^{k+j}x\|^2
$$
for every ${n\ge1}$ and every ${j\ge0}.$ Thus, since
$\|{Q_n}^\frac{1}{2}T^jx\|^2\to\|Q^\frac{1}{2}T^jx\|^2$ by (a), then
$$
\smallfrac{1}{n}{\sum}_{k=0}^{n-1}\|T^{k+j}x\|^2\to\|Q^\frac{3}{2}x\|^2,
\quad\;\hbox{for every $x\in\X$ and every $j\ge0$}.
$$
\vskip-4pt

\vskip6pt\noi
(d)
$\!$According to (c), $\|Q^\frac{1}{2}x\|^2
\kern-2pt=\kern-1pt\|Q^\frac{1}{2}T^kx\|^2
\kern-2pt\le\kern-1pt\|Q\|\kern1pt\|T^kx\|^2\kern-1pt$
for any ${k\kern-1pt\ge\kern-1pt1}$ and
$\|Q^\frac{1}{2}x\|^2\kern-2pt=\kern-1pt$
$\smallfrac{1}{n}{\sum}_{k=0}^{n-1}\|Q^\frac{1}{2}x\|^2
\kern-2pt\le\kern-1pt\|Q\|\kern1pt\smallfrac{1}{n}{\sum}_{k=0}^{n-1}\|T^kx\|^2
\kern-2pt\to\kern-1pt\|Q\|\kern1pt\|Q^\frac{1}{2}x\|^2\kern-2pt\ne\kern-1pt0$
for every ${\kern-1ptx\kern-1pt\in\kern-1pt\X\\\N(Q)}.$ So
$$
1\le\|Q\|
\quad\;\hbox{whenever}\;\quad
Q\kern-1pt\ne O.
$$
Since $Q=T^*Q\,T$, then $\|Q\|\le\|Q\|T\|^2.$ Hence ${Q\ne O}$ implies
${1\le\|T\|}$.

\vskip6pt\noi
(e)
${Q\kern-1pt=\kern-1ptO}$ trivially implies
${Q\kern1ptT\kern-1pt=\kern-1ptTQ\kern-1pt=\kern-1ptO}$, and
${Q\kern1ptT\kern-1pt=\kern-1ptO}$ implies ${Q\kern-1pt=\kern-1ptO}$ by (c)$.$
If\/ ${TQ\kern-1pt=\kern-1ptO}$, then
${0=\|{Q_n}^\frac{1}{2}TQx\|\to\|Q^\frac{3}{2}x\|}$ for every
${x\kern-1pt\in\kern-1pt\X}$ by (c) again, and so ${Q\kern-1pt=\kern-1ptO}$.

\vskip6pt\noi
(f)
Since ${O\le Q}$, then ${Q\kern1ptT=TQ}$ if and only if
${Q^\frac{1}{2}\kern1ptT=TQ^\frac{1}{2}}.$ If
${Q^\frac{1}{2}\kern1ptT=TQ^\frac{1}{2}}$, then according to (c)
it follows that
\begin{eqnarray*}
\<Qx\,;x\>
=
\|Q^\frac{1}{2}x\|^2\!
=\smallfrac{1}{n}{\sum}_{k=1}^{n-1}\|Q^\frac{1}{2}x\|^2\!
&\kern-6pt=\kern-6pt&
\smallfrac{1}{n}{\sum}_{k=1}^{n-1}\|Q^\frac{1}{2}T^kx\|^2                 \\
&\kern-6pt=\kern-6pt&
\smallfrac{1}{n}{\sum}_{k=1}^{n-1}\|T^kQ^\frac{1}{2}x\|^2
\to\|Qx\|^2
=\<Q^2x\,;x\>
\end{eqnarray*}
for every ${x\in\X}.$ So ${Q=Q^2}.$ Conversely,
${\<Q\kern1ptT^nx\,;T^nQx\>}
={\|Q^\frac{1}{2}T^nx\|^2}
={\|Q^\frac{1}{2}x\|^2}$
according to (c)$.$ Since ${Q=Q^2}$ if and only if ${Q^\frac{1}{2}=Q}$, then
we get in this case
\begin{eqnarray*}
\|(Q\kern1ptT^n-T^nQ)x\|^2
&\kern-6pt=\kern-6pt&
\|Q\kern1ptT^nx\|^2+\|T^nQx\|^2\!-2\|Q^\frac{1}{2}x\|^2                \\
&\kern-6pt=\kern-6pt&
\|T^nQx\|^2-\|Qx\|^2
\le({\sup}_n\|T^n\|^2-1)\|Qx\|^2
\end{eqnarray*}
for all $n$ and every ${x\in\X}.$ In particular, for ${n=1}$ we get for every
${x\in\X}$
$$
\|(Q\kern1ptT-TQ)x\|^2\le(\|T\|^2-1)\|Qx\|^2.
$$
Thus, asymptotically (with the assumption ${Q=Q^2}$ still in force), we get
by (c)
$$
{\lim}_n\|(Q\kern1ptT^n-T^nQ)x\|^2
={\lim}_n\|T^nQ^\frac{1}{2}x\|^2\!-\kern-1pt\|Q^\frac{1}{2}x\|^2
=0
$$
for every ${x\in\X}$ (and so $\{QT^n\kern-1pt-\kern-1ptT^nQ\}$ is a bounded
the sequence of operators disregarding whether $T$ is power bounded or not)$.$
Again, as in the proof of Theorem 6.1, ${O\le Q=Q^2\ne O}$ implies
$Q$ is an nonzero orthogonal projection, and so ${\|Q\|=1}$.

\vskip6pt\noi
(g)
The inequality is readily verified and the limit comes from (c)$:$ for any
${x\in\X}$,
\begin{eqnarray*}
\|(I-Q)\kern1pt{Q_n}^\frac{1}{2}x\|^2
&\kern-6pt=\kern-6pt&
\|{Q_n}^\frac{1}{2}x\|^2
+\|Q\kern1pt{Q_n}^\frac{1}{2}x\|^2
-2\|Q^\frac{1}{2}{Q_n}^\frac{1}{2}x\|^2                             \\
&\kern-6pt\le\kern-6pt&
\|{Q_n}^\frac{1}{2}x\|^2
+\|Q\|^2\|{Q_n}^\frac{1}{2}x\|^2
-2\|Q\|\,\|{Q_n}^\frac{1}{2}x\|^2                                   \\
&\kern-6pt=\kern-6pt&
(\|Q\|-1)^2\|{Q_n}^\frac{1}{2}x\|^2
\to(\|Q\|-1)^2\|Q^\frac{1}{2}x\|^2,
\end{eqnarray*}
and since ${I-Q}={(I+Q^\frac{1}{2})}{(I-Q^\frac{1}{2})}$ and
${I+Q^\frac{1}{2}}$ is invertible, then we get the second inequality form
the above one.

\vskip6pt\noi
(h)
As we saw above,
$\|Q\kern1pt{Q_n}^\frac{1}{2}x\|^2
=\|(I-Q)\kern1pt{Q_n}^\frac{1}{2}x\|^2
+2\|Q^\frac{1}{2}{Q_n}^\frac{1}{2}x\|^2-\|{Q_n}^\frac{1}{2}x\|^2$,
but
$2\|Q^\frac{1}{2}{Q_n}^\frac{1}{2}x\|^2-\|{Q_n}^\frac{1}{2}x\|^2
=\|(2Q-I)^\frac{1}{2}{Q_n}^\frac{1}{2}x\|^2$,
for every ${x\in\X}.$ So we get (h).

\vskip6pt\noi
(i)
As in the proof of Proposition 7.1, if ${\|T^nx\|\to0}$, then
${\frac{1}{n}{\sum}_{k=0}^{n-1}\|T^kx\|\to0}$, which means
${\|Q^\frac{1}{2}x\|=0}$ by (c) or, equivalently,
${x\in\N(Q^\frac{1}{2}})$ (i.e., ${x\in\N(Q)}\kern1pt).$ The converse requires
power boundedness$.$ If ${\frac{1}{n}{\sum}_{k=0}^{n-1}\|T^kx\|\to0}$ (i.e.,
if ${x\in\N(Q)}\kern1pt)$, then
$$
0\le{\liminf}_n\|T^nx\|
\le{\lim}_n{\inf}_j\smallfrac{1}{n}{\sum}_{k=0}^{n-1}\|T^{k+j}x\|
\le{\lim}_n\smallfrac{1}{n}{\sum}_{k=0}^{n-1}\|T^kx\|
=0
$$
(as we saw in Section 6)$.$ But if $T$ power bounded, then
${{\liminf}_n\|T^nx\|=0}$ implies ${{\lim}_n\|T^nx\|=0}$ (as we saw in the
proof of Theorem 6.1(i)\kern.5pt).

\vskip6pt\noi
(j)
Suppose ${\sup_n\|T^n\|\le\beta}.$ Again, as in the proof of Proposition 7.1,
$$
\|T^nx\|\to\beta\|x\|
\;\limply\;
\smallfrac{1}{n}{\sum}_{k=0}^{n-1}\|T^kx\|^2\to\beta^2\|x\|.
$$
According to (c),
$$
{\lim}_n\smallfrac{1}{n}{\sum}_{k=0}^{n-1}\|T^kx\|^2=\beta^2\|x\|
\iff
\|Q^\frac{1}{2}x\|^2=\beta^2\|x\|^2
\iff
\<(Q-\beta^2I)x\,;x\>,
$$
and according to (b),
$$
\<(Q-\beta^2I)x\,;x\>
\iff
\|(Q-\beta^2I)^\frac{1}{2}x\|=0
\iff
x\in\N(Q-\beta^2I)
$$
since ${\N(Q-\beta^2I)^\frac{1}{2\kern-1pt}=\N(Q-\beta^2I)}.$ Conversely,
$$
{\lim}_n\smallfrac{1}{n}{\sum}_{k=0}^{n-1}\|T^kx\|^2=\beta^2\|x\|
\;\limply\;
{\limsup}_n\|T^nx\|=\beta\|x\|.
$$
Indeed, since ${\sup_n\|T^n\|\le\beta}$, then as we saw in Section 6
$$
\beta^2\|x\|^2
={\lim}_n\smallfrac{1}{n}{\sum}_{k=0}^{n-1}\|T^kx\|^2
\le{\limsup}_n\|T^nx\|^2
\le{\sup}_n\|T^nx\|^2
\le\beta^2\|x\|^2.
$$
Thus, as in the proof of Theorem 6.1(j), for ${x\ne0}$ and since ${\beta>0}$,
$$
{\limsup}_n\|T^nx\|=\beta\|x\|
\;\limply\;
\beta\ge1
\;\;\hbox{and}\:\;
\|x\|\le{\liminf}_n\|T^nx\|.
$$
If ${\lim_n\|T^nx\|=\beta\|x\|}$ for every
${x\in\X}$, then
${\lim}_n\smallfrac{1}{n}{\sum}_{k=0}^{n-1}\|T^kx\|^2\kern-1pt=\beta^2\|x\|$
for every ${x\in\X}$ as we saw above, meaning ${Q=\beta^2I}$, which implies
${T^*T=I}$
by (c), and so ${Q_n\kern-1pt=I=Q}.$ Conversely,
if ${Q=I}$, then ${\|T^nx\|\kern-1pt=\kern-1pt\|x\|}$
for every ${x\in\X}$ by (c) again (i.e., $T$ is an
isometry), and so ${\lim_n\|T^nx\|=\beta\|x\|}$ for every ${x\in\X}$ with
${\beta=1}$.

\vskip6pt
Suppose $T$ is power bounded$.$ If $T$ is of class $C_{1{\textstyle\cdot}}$
(i.e.,if ${\|T^nx\|\not\to0}$ if ${x\ne0}$), then as we saw in the proof
of Theorem 6.1(i) ${0<\liminf_n\|T^nx\|}$ for ${x\ne0}.$ The converse is
trivial$.$ Since
$\liminf_n\|T^nx\|^2\le{\lim}_n\smallfrac{1}{n}{\sum}_{k=0}^{n-1}\|T^kx\|^2$
(as we saw in Section 6) and
${\lim}_n\smallfrac{1}{n}{\sum}_{k=0}^{n-1}\|T^kx\|^2\kern-1pt
=\|Q_nx\|^2\kern-1pt\to\|Qx\|^2\kern-1pt=\|x\|_Q^2$
by (c), then ${0<\liminf_n\|T^nx\|}$ implies ${0<\|x\|_Q}$ for ${x\ne0}$,
which means ${Q>O}.$ Thus if $T$ is power bounded, then
$$
\centerline{
${\<\,\cdot\,;\cdot\,\>_Q}$ is an inner product $\iff$ $T$ is of class
$C_{1{\textstyle\cdot}}\!\iff Q$\/ is positive.}
$$

\vskip0pt\noi
(k)
This follows as in the proof of Theorem 6.1(k) with $A_\vphi$ replaced by $Q$.

\vskip6pt
Moreover, the equivalences among the assertions (1) to (4), depend on the new
inner product ${\<\cdot\,;\cdot\>_Q}$ generated by the positive $Q$, and so
they follow by Propositions 4.1 and 4.2 by using the same argument of
Theorem 6.1, with $\vphi(\{\|T^nx\|^2\})=\|x\|_\vphi^2$ replaced by
$\lim_n\|{Q_n}^\frac{1}{2}x\|^2
={\lim}_n\smallfrac{1}{n}{\sum}_{k=0}^{n-1}\|T^kx\|^2
\kern-1pt\to\|Qx\|^2\kern-1pt=\|x\|_Q^2$.
\end{proof}
\vskip-4pt

\vskip3pt\noi
\begin{remark}
Even a power unbounded operator may have a Ces\`aro asymptotic limit (see,
e.g., \cite[Example 3]{Geh2}), while there is no $\vphi$-asymptotic limit for
power unbounded operators$.$ From now on suppose $T$ is power bounded$.$

\vskip4pt\noi
(a)
Thus for every Banach limit $\vphi$ there exists a $\vphi$-asymptotic limit
$A_\vphi$ for $T.$ Even in this case of a power bounded operator, the
Ces\`aro asymptotic limit $Q$ may not exist (even in the weak sense; see
\cite[Example 2]{Geh2}).

\vskip4pt\noi
(b)
Moreover, even when $Q$ exists it may not coincide with $A_\vphi.$ Indeed, it
was exhibited in \cite[Example 1]{Geh2} a power bounded unilateral weighted
shift $T$ such that ${\|T^n\|=\beta=\sqrt2}$ for all $n$ and ${\|T^ne_1\|^2}$
is either $\beta^2$ or $1$ depending on $n$, with Ces\`aro asymptotic limit
${Q=I}$, which does not coincide with an arbitrary $\vphi$-asymptotic limit
$A_\vphi$ (i.e., ${Q\ne A_\vphi}$ for a specific Banach limit $\vphi$ ---
actually, there exist Banach limits $\vphi$ for which
${\|{A_\vphi}^{_{\!\scriptstyle\frac{1}{2}}}e_1\|^2}$ lies anywhere in the
interval ${[1,\beta^2]})$.

\vskip4pt\noi
(c)
As we have seen in Theorems 6.1(d) and 7.1(d), if the asymptotic limits are
not null, then for every Banach limit $\vphi$ we get
${1\kern-1pt\le\kern-1pt\|A_\vphi\|}$, $\,{1\kern-1pt\le\kern-1pt\|Q\|}$, and
${1\kern-1pt\le\kern-1pt\|T\|}.$ These norms, however, are not related$.$
For instance, if ${T={\rm shift}\{\beta,1,1,1,\dots\}}$ is the unilateral
weighted shift with ${\beta>1}$ as in the proofs of Theorems 6.1(b) and
7.1(b), then ${\|A_\vphi\|\kern-1pt=\kern-1pt\|Q\|\kern-1pt=\kern-1pt\beta^2}$
and ${\|T\|\kern-1pt=\kern-1pt\beta}.$ On the other hand, if
$T=\big(\smallmatrix{0 & \beta \cr
                     0 & 0     \cr}\big)\oplus I$,
then $T^n\kern-1pt={O\oplus I}$ for all ${n\ge2}$ and
${A_\vphi\kern-1pt=Q=O\oplus I}$ for all Banach limits $\vphi$ and so
${\|T\|=\beta}$ with
${\|A_\vphi\|\kern-1pt=\kern-1pt\|Q\|\kern-1pt=\kern-1pt1}.$ Actually,
as we saw in item (b) above, it was exhib\-ited in \cite[Example 1]{Geh2} a
power bounded unilateral weighted shift $T$ such that there is a maximum
Banach limit $\vphi_+$ for which ${\|A_{\vphi_+}\|\kern-1pt\ge\kern-1pt2}$,
while ${\|T\|\kern-1pt=\kern-1pt{\sqrt 2}}$ and ${\|Q\|\kern-1pt=\kern-1pt1}$.

\vskip4pt\noi
(d) 
The inclusions in Theorems 6.1(j) and 7.1(j) may also be proper (e.g.,
for the unilateral weighted shift $T$ form \cite[Example 1]{Geh2}, as in
item (b) above, ${\beta^2=2}$ and ${Q=I}$ so that ${\N(\beta^2I-Q)=\0}$ while
$\|T^ne_1\|$ oscillates between $1$ and $\beta$).
\end{remark}
\vskip-2pt

\vskip6pt
For a power bounded operator on a finite-dimensional space, the Ces\`aro
asymptotic limit $Q$ exists and coincides with the $\vphi$-asymptotic limit
$A_\vphi$ for every Banach limit $\vphi$ \cite[Theorem 2.1]{Geh2}$.$ The next
theorem gives a condition for ${Q=A_\vphi}$ on an infinite-dimensional
space$.$ As we saw in the proof of Theorem 7.1(c), if the sequence $\{Q_n\}$
of Ces\`aro means converges weakly to, say $Q$, then the sequence
$\{T^{*j}Q_nT^j\}$ of Ces\`aro means converges weakly (again to $Q$) for every
positive integer $j.$ If such weak convergence holds uniformly in $j$, then
${Q=A_\vphi}$ for all Banach limits $\vphi$.

\vskip3pt\noi
\begin{theorem}
If\/ $T$ is a Hilbert-space power bounded operator for which the sequence\/
$\{T^{*j}Q_nT^j\}$ of Ces\`aro means converges weakly and uniformly in\/ $j$,
$$
T^{*j}Q_nT^j=\smallfrac{1}{n}{\sum}_{k=0}^{n-1}T^{*k+j}T^{k+j}\wconv Q
\quad\;\hbox{uniformly in\/ $j$},
$$
then the Ces\`aro asymptotic limit coincides with the\/ $\vphi$-asymptotic
limit,
$$
Q=A_\vphi,
$$
\vskip-4pt\noi
for all Banach limits\/ ${\vphi\!:\ell_+^\infty\!\to\CC}$.
\end{theorem}

\vskip0pt\noi
\begin{proof}
By Theorem 7.1(c), ${Q_n\wconv Q}$ if and only if ${T^{*j}Q_nT^j\wconv Q}$
which means
$$
\smallfrac{1}{n}{\sum}_{k=0}^{n-1}\|T^{k+j}x\|^2\to\|Q^\frac{1}{2}x\|^2
\quad\;\hbox{for every\/ ${x\in\X}$ and every $j\ge0$},
$$
as ${Q\ge O}.$ If the weak convergence of $\{T^{*j}Q_nT^j\}$ holds uniformly
in $j$, then so does the above convergence$.$ But the real-valued sequence
$\big\{\smallfrac{1}{n}{\sum}_{k=0}^{n-1}\|T^{k+j}x\|^2\big\}$ of Ces\`aro
means converges uniformly in $j$ if and only if all Banach limits
${\vphi\in\kern-1pt\X^*}\kern-1pt$ coincide at the sequence
$\kern-1pt\{\|T^nx\|^2\kern-1pt\}\kern-1pt$ and are equal to
$\|Q^\frac{1}{2}x\|^2\kern-1pt$ \cite[\!Theorem 1]{Lor} (also \cite{Such});
that is,
$$
\vphi(\{\|T^nx\|^2\})=\|Q^\frac{1}{2}x\|^2
$$
for all Banach limits ${\vphi\!:\ell_+^\infty\!\to\CC}.$ In particular, this
holds for the arbitrary Banach limit $\vphi$ of Theorem 6.1 since $T$ is
power bounded$.$ For that Banach limit we got
$$
\vphi(\{\|T^nx\|^2)=\|{A_\vphi}^{_{\!\scriptstyle\frac{1}{2}}}x\|^2
$$
where ${A_\vphi\kern-1pt\ge O}$ is the $\vphi$-asymptotic limit of $T$
(associated with $\vphi).$ Thus
$\|Q^\frac{1}{2}x\|^2=\|{A_\vphi}^{_{\!\scriptstyle\frac{1}{2}}}x\|^2$
or, equivalently, ${\<(Q-A_\vphi)x,;x\>=0}$, for every ${x\in\X}.$ This means
$$
Q=A_\vphi
$$
(either because the Hilbert space is complex or because ${Q-A_\vphi}$ is
self-adjoint).
\end{proof}
\vskip-2pt

\vskip3pt\noi
\begin{remark}
(a) A class of operators that satisfies the assumption of Theorem 7.2 is the
class of quasinormal operators$.$ A Hilbert-space operator $T$ is quasinormal
if it commutes with ${T^*T}.$ If ${T\kern-1pt\in\kern-1pt\BX}$ is quasinormal
on a Hilbert space $\X$, then by two trivial inductions we get
${T^*T\,T^k\kern-1pt=T^kT^*T}$ for every ${k\ge1}$, and consequently
${T^{*k}T^k}\kern-1pt={(T^*T)^k}$ for every ${k\ge1}.$ This in fact is
equivalent to quasinormality --- see, e.g., \cite[Proposition 13]{Jib} and
\cite[Theorem 3.6]{JJS}$.$ Therefore
\vskip4pt\noi
{\narrower
\it there is an operator $S$ for which ${T^{*k}T^k}\!={S^k}\kern-1pt$ for
every ${k\kern-1pt\ge\kern-1pt1}$ if and only if $\kern1ptT$ is quasinormal,
and such an operator is unique and given by
${S\kern-1pt=|T|^2\!=\kern-1ptT^*T}$.
\vskip4pt}\noi
In this case, ${T^{*k+j}T^{k+j}\kern-1pt=(T^*T)^{k+j}\kern-1pt=S^{k+j}}$ for
every nonnegative integers $j,k.$ If $T$ is power bounded, then so is $S$, and
the Mean Ergodic Theorem for power bounded operators (which holds in reflexive
Banach spaces --- see, e.g., \cite[Corollary VIII.5.4]{LO1}) ensures strong
convergence for the sequence of Ces\`aro means
$\big\{\kern-1pt\smallfrac{1}{n}{\sum}_{k=0}^{n-1}\kern1ptS^k\kern-1pt\big\}$
whose strong limit $Q$ lies in $\BX$ by the Banach--Steinhaus Theorem$.$ Thus
$$
Q_n=\smallfrac{1}{n}{\sum}_{k=0}^{n-1}T^{*k}T^k
=\smallfrac{1}{n}{\sum}_{k=0}^{n-1}\kern1ptS^k\sconv Q,
$$
where ${Q\ge O}$ is the Ces\`aro asymptotic limit of $T$ (cf. proof
Theorem 7.1)$.$ Hence
$$
\smallfrac{1}{n}{\sum}_{k=0}^{n-1}S^{k+j}
=\smallfrac{1}{n}{\sum}_{k=0}^{n-1}T^{*k+j}T^{k+j}
=T^{*j}Q_nT^j
\sconv T^{*j}QT^j=Q=S^jQ=QS^j
$$
for every $j$ according to Theorem 7.1 (as strong convergence implies
weak conver\-gence)$.$ Take an arbitrary ${x\in\X}.$ By the above strong
convergence
\begin{eqnarray*}
{\sup}_j\Big\|
\Big(\smallfrac{1}{n}{\sum}_{k=0}^{n-1}T^{*k+j}T^{k+j}-Q\Big)x\Big\|
&\kern-6pt=\kern-6pt&
{\sup}_j\Big\|S^j\Big(\smallfrac{1}{n}{\sum}_{k=0}^{n-1}S^k-Q\Big)x\Big\|  \\
&\kern-6pt\le\kern-6pt&
{\sup}_j\|S^j\|\,
\Big\|\Big(\smallfrac{1}{n}{\sum}_{k=0}^{n-1}S^k-Q\Big)x\Big\|\to0.
\end{eqnarray*}
\vskip-2pt\noi
Then $\smallfrac{1}{n}{\sum}_{k=0}^{n-1}T^{*k+j}T^{k+j}\sconv Q$ uniformly in
$j.$ Hence (again, since strong convergence implies weak convergence)
$\smallfrac{1}{n}{\sum}_{k=0}^{n-1}T^{*k+j}T^{k+j}\wconv Q$ uniformly in $j.$
(Indeed,
$\sup_j\big|\big\<\big(\smallfrac{1}{n}
{\sum}_{k=0}^{n-1}T^{*k+j}T^{k+j}-Q\big)x\,;x\big>\big|
\le\sup_j\big\|
\big(\smallfrac{1}{n}{\sum}_{k=0}^{n-1}T^{*k+j}T^{k+j}
-Q\big)x\big\|\kern1pt\|x\|$.)
So
$$
\smallfrac{1}{n}{\sum}_{k=0}^{n-1}\|T^{k+j}x\|^2\to\|Q^\frac{1}{2}x\|^2
\quad\;\hbox{uniformly\/ in $j$}.
$$
Thus, according to Theorem 7.2, ${Q=A_\vphi}$ for all Banach limits $\vphi$,
where $A_\vphi$ is the $\vphi$-asymptotic limit for the power bounded
operator $T$ as in Theorem 6.1.

\vskip6pt\noi
(b) A normed-space operator $T$ is normaloid if ${\|T^n\|=\|T\|^n}$ for every
integer $n\ge0.$ By the Gelfand--Beurling formula, on a complex Banach-space
a normaloid is an operator $T$ for which spectral radius coincides with norm:
$r(T)=\|T\|.$ Since power boundedness implies $r(T)\le1$, then it follows at
once that
\vskip4pt\noi
{\narrower\narrower
\it a power bounded operator is normaloid if and only if it is a normaloid
contraction\/.
\vskip4pt}\noi
(In fact, {\it if a normaloid operator is similar to a power bounded operator,
then it is a contraction}\/ \cite[Proposition 1]{Kub2}.) Quasinormal is a
class of operators including the normal operators and the isometries, and it
is included in the class of subnormal operators, which is included in the
class of hyponormal operators, which in turn is included in the class of
paranormal operators, which are all normaloid$.$ So all these Hilbert-space
normaloid operators, when power bounded, are contractions and so they
naturally fit to Proposition 7.1 (and consequently they trivially fit to
Theorem 7.2 --- see Remark 6.2).
\end{remark}
\vskip-2pt

\vskip3pt\noi
\begin{corollary}
Let\/ $T$ be a Hilbert-space power bounded operator. If the sequence\/
$\{Q_n\}$ of Ces\`aro means converges uniformly,
$$
Q_n=\smallfrac{1}{n}{\sum}_{k=0}^{n-1}T^{*k}T^k\uconv Q,
$$
then the Ces\`aro asymptotic limit coincides with the\/ $\vphi$-asymptotic
limit,
$$
Q=A_\vphi,
$$
\vskip-4pt\noi
for all Banach limits\/ ${\vphi\!:\ell_+^\infty\!\to\CC}$.
\end{corollary}

\vskip0pt\noi
\begin{proof}
Consider the setup of Theorem 7.1$.$ Recall that ${Q=T^{*j}QT^j}$ for
every ${j\ge1}$. If
${Q_n\!=\kern-1pt\smallfrac{1}{n}{\sum}_{k=0}^{n-1}T^{*k}T^k\!\uconv Q}$,
then
$$
{\sup}_j\Big\|\smallfrac{1}{n}{\sum}_{k=0}^{n-1}T^{*k+j}T^{k+j}-Q\Big\|
\le
{\sup}_j\|T^j\|^2\Big\|\smallfrac{1}{n}{\sum}_{k=0}^{n-1}T^{*k}T^{k}-Q\Big\|
\to0.
$$
Thus $\smallfrac{1}{n}{\sum}_{k=0}^{n-1}T^{*k+j}T^{k+j}\uconv Q$ (and so
$\smallfrac{1}{n}{\sum}_{k=0}^{n-1}T^{*k+j}T^{k+j}\wconv Q$) uniformly in $j.$
If $T$ is power bounded, then ${Q=A_\vphi}$ for all $\vphi$ by Theorem 7.2$.$
\end{proof}
\vskip-2pt

\vskip6pt
For instance, let $T$ is a uniformly stable noncontraction (i.e.,
${r(T)<1<\|T\|}$) acting on any Hilbert space $\X.$ Then ${T^n\!\uconv O}$
(so that $T$ is power bounded) or, equivalently, ${T^{*n}T^n\!\uconv O}$,
and so $Q_n=\smallfrac{1}{n}{\sum}_{k=0}^{n-1}T^{*k}T^k\!\uconv Q=O=A_\vphi$
for all Banach \hbox{limits}\/ ${\vphi\!:\ell_+^\infty\!\to\CC}$ (in
accordance with Corollary 7.1).

\vskip3pt\noi
\begin{remark}
If $T$ is a power bounded operator on a finite-dimensional space, then
${Q=A_\vphi}$ for all Banach limits\/ ${\vphi\!:\ell_+^\infty\!\to\CC}.$
Indeed, for power bounded operators on a finite-dimensional space (where
weak, strong, and uniform convergences coincide), the Ces\`aro asymptotic
limit exists \cite[Theorem 2.1]{Geh2}$.$ Thus ${Q=A_\vphi}$ for all Banach
limits $\vphi$ by Corollary 7.1.
\end{remark}

\section*{Acknowledgment}

We thank Gy\"orgy Geh\'er for enlightening discussions on Banach limits.

\bibliographystyle{amsplain}

\end{document}